# BAYESIAN-MOTIVATED TESTS OF FUNCTION FIT AND THEIR ASYMPTOTIC FREQUENTIST PROPERTIES

By Marc Aerts[1], Gerda Claeskens[2] and Jeffrey D. Hart[3]

*Limburgs Universitair Centrum, K. U. Leuven and Texas A&M University*

We propose and analyze nonparametric tests of the null hypothesis that a function belongs to a specified parametric family. The tests are based on BIC approximations, $\pi_{\mathrm{BIC}}$, to the posterior probability of the null model, and may be carried out in either Bayesian or frequentist fashion. We obtain results on the asymptotic distribution of $\pi_{\mathrm{BIC}}$ under both the null hypothesis and local alternatives. One version of $\pi_{\mathrm{BIC}}$, call it $\pi_{\mathrm{BIC}}^*$, uses a class of models that are orthogonal to each other and growing in number without bound as sample size, $n$, tends to infinity. We show that $\sqrt{n}(1 - \pi_{\mathrm{BIC}}^*)$ converges in distribution to a stable law under the null hypothesis. We also show that $\pi_{\mathrm{BIC}}^*$ can detect local alternatives converging to the null at the rate $\sqrt{\log n/n}$. A particularly interesting finding is that the power of the $\pi_{\mathrm{BIC}}^*$-based test is asymptotically equal to that of a test based on the maximum of alternative log-likelihoods.

Simulation results and an example involving variable star data illustrate desirable features of the proposed tests.

**1. Introduction.** Consider a model in which the observed data vector **Y** has distribution $f(\mathbf{y}, g, \boldsymbol{\eta})$, where $f$ is known, $g$ is an unknown function and $\boldsymbol{\eta}$ is a vector of unknown nuisance parameters. We wish to test the null hypothesis that $g$ is in a specified parametric family $\mathcal{G} = \{g(\cdot; \boldsymbol{\theta}) : \boldsymbol{\theta} \in \boldsymbol{\Theta}\}$ against the nonparametric alternative that $g \notin \mathcal{G}$. This paper proposes a Bayes-inspired test of such a hypothesis. A version of the test was proposed by Hart (1997) in the special case of checking the fit of a parametric regression model. The idea is simple. Consider a sequence of models for $g$ of

Received July 2003; revised April 2004.

[1]Supported by the Belgian IAP network nr P5/24 of the Belgian Government (Belgian Federal Science Policy Office).

[2]Supported in part by NSF Grant DMS-02-03884 and by the Belgian Federal Science Policy Office.

[3]Supported in part by NSF Grant DMS-99-71755.

varying dimensions, one of which is the parametric (or null) model whose fit is to be tested. The posterior probability, $\pi_n$, of the null model is computed, and if this probability is sufficiently low, the null model is rejected. This test may be carried out in either Bayesian or frequentist fashion. One may determine a sequence of constants $a_n$ such that $a_n(1 - \pi_n)$ converges in distribution to a nondegenerate random variable when $H_0$ is true and the sample size $n$ tends to $\infty$. This allows the frequentist to conduct a valid large sample test of given size based on $a_n(1 - \pi_n)$. On the other hand, a Bayesian may simply wish to make a decision based on the value of $\pi_n$, irrespective of an a priori type I error probability.

The idea of using a Bayesian-motivated statistic in frequentist fashion is not new. Good (1957) proposed that the distribution of a Bayes factor be computed on the assumption that a sharp null hypothesis is true, and $P$-values corresponding to the Bayes factor be used as a significance criterion. Good (1992) gives an extensive review of compromises between Bayesian and non-Bayesian methodologies.

Lack-of-fit and goodness-of-fit tests based on orthogonal series expansions and/or smoothing ideas have received considerable attention in the last fifteen or so years. Many references to this work may be found in the book of Hart (1997). Seminal references on series-based goodness-of-fit tests, that is, so-called smooth tests, are Neyman (1937) and Rayner and Best (1989, 1990). More recently, Ledwina (1994) and Fan (1996) have proposed adaptive versions of Neyman's smooth test. Eubank and Hart (1992) and Aerts, Claeskens and Hart (1999) have studied the so-called order selection test in the contexts of regression and general likelihood models, respectively. A nonparametric Bayesian goodness-of-fit test has been proposed by Verdinelli and Wasserman (1998).

The rest of the paper is organized as follows. Section 2 considers frequentist and formal Bayesian versions of the proposed test, and discusses the choice of alternative models and specification of priors. Section 3.1 summarizes a simulation study comparing the power of our test with other omnibus lack-of-fit tests. In Section 3.2 our methods are applied to the problem of testing for a trend in the sequence of times between maximum brightnesses of the long-period variable star Omicron Ceti or Mira. Section 4 presents our theoretical results on the asymptotic frequentist properties of the proposed tests. Finally, the Appendix contains mathematical details and proofs of the theorems.

**2. Test procedures.** To reiterate, we assume that observed data $\mathbf{Y}$ have distribution $f(\mathbf{y}, g, \boldsymbol{\eta})$ for some function $g$ and vector of parameters $\boldsymbol{\eta}$. We wish to test the hypothesis, call it $H_0$, that the function $g$ lies in the parametric family of functions $\mathcal{G}$. The model which assumes that $H_0$ is true will be called $M_0$. We consider a collection of alternative models denoted $M_1, \ldots, M_K$, where each $M_i$ corresponds to a different parametric specification for the function $g$. These models need not be nested within each other. Since we wish our test of $H_0$ to be nonparametric, $K$ should be fairly "large" and the union of $M_0, M_1, \ldots, M_K$ should come close to spanning the space of all possibilities for $g$. Indeed, we can envision $K$ growing with the number of observations in $\mathbf{Y}$ in such a way that, asymptotically, the models under consideration *do* span all the possibilities.

Our tests of $H_0$ are based on a posterior probability for $M_0$ or on an approximation to that probability. These tests run the gamut from a purely Bayesian approach based on informative priors to a purely frequentist one that involves no prior specification at all. In any case, our tests take the form

"reject $H_0$ when $\pi_n \stackrel{\mathrm{def}}{=} P(M_0|\mathbf{y})$ is sufficiently small."

A Bayesian will make a decision, or perhaps abstain from doing so, by simply examining $\pi_n$ and/or a Bayes factor. On the other hand, a frequentist will wish to determine the sampling distribution of $P(M_0|\mathbf{Y})$ on the assumption that $H_0$ is true, and then reject $H_0$ at level of significance $\alpha$ if and only if $\pi_n$ is smaller than an $\alpha$ quantile of this distribution. The frequentist may well regard $\pi_n$ differently than a Bayesian. The latter views $\pi_n$ as the probability that $H_0$ is true in light of the observed data, whereas the former may simply view it as a statistic that contains evidence about the hypotheses of interest.

In Section 2.1 we turn to the question of choosing alternative models $M_1, M_2, \ldots$, a question of relevance to both Bayesians and frequentists. Section 2.2 considers a formal Bayesian version of the proposed test, including a discussion of noninformative priors for the models $M_0, M_1, \ldots$. An asymptotic version of the test requiring no specification of priors is introduced in Section 2.3.

2.1. *Alternative models.* We shall consider two main types of alternative models: those which are guaranteed to contain the true function $g$ (at least in a limiting sense) and those which do not necessarily contain $g$ but nonetheless lead to a consistent test for virtually any $g$. An example will be helpful to illustrate these two types. In the sequel, the model for $g$ corresponding to probability model $M_j$ will be denoted $g_j$. Suppose that $g$ and each member of $\mathcal{G}$ are continuous functions over the interval $[0, 1]$, which means that $g$ can be written as

$$g(x) = g(x; \boldsymbol{\theta}) + \delta(x),$$

where $\delta$ has the Fourier series representation

$$\delta(x) = \sum_{j=0}^{\infty} \alpha_j \cos(\pi j x), \qquad 0 \leq x \leq 1,$$

for constants $\alpha_0, \alpha_1, \ldots$. This representation for $g$ suggests that we take

$$(1) \qquad g_j(x; \boldsymbol{\theta}, \boldsymbol{\alpha}_j) = g(x; \boldsymbol{\theta}) + \sum_{k=0}^{j} \alpha_k \cos(\pi k x),$$

where $\boldsymbol{\alpha}_j = (\alpha_0, \ldots, \alpha_j)^T$. Of course, this model could be modified to suit a given situation. For example, if $g$ is a regression function and the model $\mathcal{G}$ contains an intercept, then the constant term $\alpha_0$ should be eliminated from $g_j$. Another model that would be useful for cases where $g$ is inherently positive is

$$g_j(x; \boldsymbol{\theta}, \boldsymbol{\alpha}_j) = g(x; \boldsymbol{\theta}) \exp\left[\sum_{k=0}^{j} \alpha_k \cos(\pi k x)\right].$$

Other basis functions can be used as well; popular examples include wavelets and orthogonal Legendre or Hermite polynomials. Wavelet packets would be particularly attractive when the most parsimonious basis is unknown to the investigator.

As $j \to \infty$, functions of the form (1) span the space of all functions that are continuous on $[0, 1]$. In many settings this property is enough to ensure that there exist tests based on the models $M_1, \ldots, M_K$ that are consistent against any continuous alternative to $H_0$, so long as $K$ tends to $\infty$ at an appropriate rate with the sample size. An example of such a test is given in Aerts, Claeskens and Hart (1999).

On the other hand, it is possible to construct consistent tests based on sequences of models that do not contain, even in the limit, the true function $g$. Such sequences can have certain advantages when using the methodology proposed in this paper. For our tests to be consistent, it is usually enough that the best approximation to $g$ among the models entertained is not in $\mathcal{G}$. Again suppose that $g$ is a function defined on $[0, 1]$. Two candidates for $g_j$ are

$$(2) \qquad g(x; \boldsymbol{\theta}) + \alpha_j \cos(\pi j x) \quad \text{and} \quad g(x; \boldsymbol{\theta}) \exp[\alpha_j \cos(\pi j x)].$$

Now, if $g$ is not in $\mathcal{G}$, but is continuous, then, generally speaking, there will exist a $k$ such that the MLE of $\alpha_k$ in $g(x; \boldsymbol{\theta}) + \alpha_k \cos(\pi k x)$ consistently estimates a nonzero quantity [White (1994)]. Such a property implies the existence of a consistent test.

The alternative models considered could be more or less arbitrary. For example, in the situation discussed immediately above we could entertain all models of the form

$$g(x; \boldsymbol{\theta}) + \sum_{k \in \mathcal{K}} \alpha_k \cos(\pi k x),$$

where $\mathcal{K}$ is a subset of $0, 1, \ldots, K$ for some $K$. If $K$ grows with sample size, such alternatives are problematic in that the number of models that must be fitted is $2^{K+1}$, which becomes prohibitively large very quickly.

In the sequel we will mainly be concerned with two classes of alternative models, ones that are *nested* and ones we shall call *singletons* that contain only one more parameter than $M_0$. Nested models are such that $M_j$ is a special case of $M_{j+1}$ for $j = 0, 1, \ldots$, while singletons contain $M_0$ but are not nested within each other.

2.2. *Formal Bayes tests.* Corresponding to model $M_j$, $j = 0, 1, \ldots, K$, are the nuisance parameters $\boldsymbol{\eta}$, parameters $\boldsymbol{\theta}$ and $\boldsymbol{\alpha}_j$ that specify $g$, and the dimension of $(\boldsymbol{\theta}, \boldsymbol{\alpha}_j, \boldsymbol{\eta})$, denoted $m_j$. The likelihood function for $M_j$ is $L(\boldsymbol{\theta}, \boldsymbol{\alpha}_j, \boldsymbol{\eta})$. Let $p_j$ be the prior probability of the $j$th model, and $\pi_j(\boldsymbol{\theta}, \boldsymbol{\alpha}_j, \boldsymbol{\eta})$ the conditional prior density of $(\boldsymbol{\theta}, \boldsymbol{\alpha}_j, \boldsymbol{\eta})$ given that the true model is $M_j$. The posterior probability of $M_0$, that is, the basis of our test of $H_0$, is

$$P(M_0|\mathbf{y}) = \frac{p_0 \int L(\boldsymbol{\theta}_0, \boldsymbol{\eta}) \pi_0(\boldsymbol{\theta}_0, \boldsymbol{\eta}) \, d\boldsymbol{\theta}_0 \, d\boldsymbol{\eta}}{\sum_{j=0}^{K} p_j \int L(\boldsymbol{\theta}, \boldsymbol{\alpha}_j, \boldsymbol{\eta}) \pi_j(\boldsymbol{\theta}, \boldsymbol{\alpha}_j, \boldsymbol{\eta}) \, d\boldsymbol{\theta} \, d\boldsymbol{\alpha}_j \, d\boldsymbol{\eta}}$$

$$= \left\{ 1 + \sum_{j=1}^{K} \frac{p_j}{p_0} \cdot B_j \right\}^{-1},$$

where $B_j$ is the Bayes factor of $M_j$ to $M_0$, that is,

$$B_j = \frac{\int L(\boldsymbol{\theta}, \boldsymbol{\alpha}_j, \boldsymbol{\eta}) \pi_j(\boldsymbol{\theta}, \boldsymbol{\alpha}_j, \boldsymbol{\eta}) \, d\boldsymbol{\theta} \, d\boldsymbol{\alpha}_j \, d\boldsymbol{\eta}}{\int L(\boldsymbol{\theta}_0, \boldsymbol{\eta}) \pi_0(\boldsymbol{\theta}_0, \boldsymbol{\eta}) \, d\boldsymbol{\theta}_0 \, d\boldsymbol{\eta}}.$$

In a subjective Bayesian analysis, the prior probabilities $p_j$ and prior distributions $\pi_j$, $j = 0, 1, \ldots, K$, are chosen to represent the investigator's degree of belief in the various models and the parameters therein. A Bayesian who wishes to do an analysis independent of his or her own prior beliefs may wish to use noninformative priors. In our setting, it is necessary to formulate such priors for the parameters in each of the models $M_0, \ldots, M_K$ and also to specify "vague" prior probabilities over these models. We have little to say here about the former problem since much has already been written about it. There has been much debate about what is the most appropriate noninformative or reference prior in a given situation, and, indeed, about whether or not *any* prior can truly express ignorance about the underlying parameters. Rather than entering this debate, we refer the interested reader to the excellent review article of Kass and Wasserman (1996) for further discussion of the problem and many relevant references.

We turn now to the question of assigning vague prior probabilities to the models $M_0, \ldots, M_K$. One possibility is to simply give each model the same probability of $1/(K + 1)$. In as much as $H_0$ has some special significance (scientific or otherwise), there may be a prevailing a priori degree of belief in it, expressed by $p_0 = \pi$. In this case we could take $p_j = (1 - \pi)/K$, $j = 1, \ldots, K$, to express lack of preference for any alternative model.

For some choices of alternative models it is debatable whether assigning them equal probabilities is really noninformative. When the models are nested with $m_0 < m_1 < \cdots$, one could argue that it is natural to put smaller prior probabilities on the models of larger dimension. Jeffreys (1961) proposed using the improper prior $p_j = 1/(j + 1)$, $j = 0, 1, \ldots$, for such problems. A proper noninformative prior for the positive integers was proposed by Rissanen (1983).

Sometimes one may consider more than one model having a given dimension. If the distinct model dimensions are $m_0 < m_1 < \cdots$, then we may assign prior probability of $2^{-M(j+1)}$ to the collection of models having dimension $m_j$ and equal probability to each individual model of that dimension. Such a scheme has been proposed by Berger and Pericchi (1996).

It is of some interest to know what form $\pi_n$ takes in various cases. Hart (1997) obtains an explicit expression for a very accurate approximation to $\pi_n$ in a regression context where one tests the hypothesis that the regression function is flat. In most cases, though, it will not be possible to write this probability as an explicit function of the data. Numerical integration or use of MCMC methods will then be needed to compute $\pi_n$.

2.3. *Tests free of prior specification.* Let $m(\mathbf{y})$ be the marginal distribution of the data $\mathbf{Y}$. In deriving the well-known BIC for selecting model dimension, Schwarz (1978) showed that in exponential family models

$$\log(P(M_j|\mathbf{y})) \approx \log L_j - \tfrac{1}{2} m_j \log n - \log(m(\mathbf{y}))$$
$$= \mathrm{BIC}_j - \log(m(\mathbf{y})),$$

where $n$ denotes the dimension of $\mathbf{y}$, $m_j$ is model dimension and $L_j$ is the likelihood function for model $j$ evaluated at the MLE. Applying this approximation to our test statistic $P(M_0|\mathbf{y})$ yields

$$P(M_0|\mathbf{y}) \approx \frac{1}{1 + \sum_{j=1}^{K} \exp(\mathrm{BIC}_j - \mathrm{BIC}_0)} \overset{\text{def}}{=} \pi_{\mathrm{BIC}}.$$

Perhaps the most interesting aspect of this approximation, especially for a frequentist, is that it is completely free of prior probabilities. The statistic $\pi_{\mathrm{BIC}}$ would seem to be attractive to frequentists and Bayesians alike. The frequentist will appreciate the fact that $\pi_{\mathrm{BIC}}$ requires no specification of priors and is thus immediately usable as a test of $H_0$ versus general alternatives. For a Bayesian, $\pi_{\mathrm{BIC}}$ can serve as a rough and ready approximation to the posterior probability of $M_0$ when the sample size is large, a property established in various contexts by Schwarz (1978), Haughton (1988), Kass and Raftery ([1995](#)) and Kass and Wasserman ([1995](#)). The reader is cautioned, however, that $\pi_{\mathrm{BIC}}$ will not always be an adequate approximation. This is especially true in small to moderate sample sizes. Furthermore, the approximation can be poor depending on the type of prior distribution used for the parameters of the models $M_0, M_1, \ldots, M_K$. For more on this last point, the reader is referred to Kass and Wasserman ([1995](#)).

2.4. *A frequentist test.* Let $\mathcal{A} = \{M_1, \ldots, M_K\}$ be a collection of models, each of which contains the null model $M_0$ as a special case. We consider the test that rejects the null hypothesis for large values of $1 - \pi_{\mathrm{BIC}}$, where

$$\pi_{\mathrm{BIC}} = \left\{ 1 + \sum_{j=1}^{K} n^{-(1/2)(m_j - m_0)} \exp\{\mathcal{L}_j / 2\} \right\}^{-1},$$

$\mathcal{L}_j$ is the log-likelihood ratio $2 \log(L_j / L_0)$, and $m_j$ denotes the number of parameters in $M_j$, $j = 0, \ldots, K$. Some of the theory to be developed later assumes that the model is of generalized linear form. In this case, the observed data are $(\mathbf{x}_1, Y_1), \ldots, (\mathbf{x}_n, Y_n)$, where each $\mathbf{x}_i$ is a vector of covariates and each $Y_i$ a scalar response. Assuming the covariates to be fixed and the observations to be independent, the log-likelihood function has the form

$$\ell(g, \eta) = \sum_{i=1}^{n} [Y_i g(\mathbf{x}_i) - b(g(\mathbf{x}_i))] / a(\eta) + c(Y_i, \eta),$$

where $a(\cdot), b(\cdot)$ and $c(\cdot)$ are known functions, $g$ is an unknown function and $\eta$ an unknown dispersion parameter; see, for example, McCullagh and Nelder ([1989](#)). We consider testing the null hypothesis

(3) $$H_0 : g(x) = \sum_{j=1}^{p} \theta_j \gamma_j(x) \stackrel{\text{def}}{=} g(x; \boldsymbol{\theta}),$$

where $\gamma_1, \ldots, \gamma_p$ are known functions and $\boldsymbol{\theta} = (\theta_1, \ldots, \theta_p)^T$ an unknown parameter vector. The asymptotic maximizer of the expected log-likelihood

$$\frac{1}{n} \sum_{i=1}^{n} [b'(g(x_i)) g(x_i; \boldsymbol{\theta}) - b(g(x_i; \boldsymbol{\theta}))]$$

with respect to $\boldsymbol{\theta}$ is denoted $\boldsymbol{\theta}_0 = (\theta_{10}, \ldots, \theta_{p0})^T$, which is the true parameter vector when $H_0$ is true and provides a best null approximation to $g$ when $H_0$ is false.

Our most general alternatives $M_j$, $j = 1, 2, \ldots,$ are of the form

$$g_j(x) = \sum_{i=1}^{p} \theta_i \gamma_i(x) + \sum_{i=1}^{j} \alpha_i v_i(x)$$

for appropriate functions $v_1, v_2, \ldots$. To produce test statistics that are meaningful and powerful, we insist that the $v_j$'s be orthonormal in the following sense:

$$(4) \qquad \sum_{i=1}^{n} v_j(\mathbf{x}_i) \gamma_k(\mathbf{x}_i) b''(g(\mathbf{x}_i; \boldsymbol{\theta}_0)) = 0, \qquad j = 1, 2, \ldots; k = 1, \ldots, p,$$

and for all $j, k \geq 1$,

$$(5) \qquad \frac{1}{n} \sum_{i=1}^{n} b''(g(\mathbf{x}_i; \boldsymbol{\theta}_0)) v_j(\mathbf{x}_i) v_k(\mathbf{x}_i) = \begin{cases} 1, & j = k, \\ 0, & j \neq k. \end{cases}$$

In practice, we may achieve an approximation to (4) and (5) by proceeding as follows. First, let $(\hat{\boldsymbol{\theta}}_0, \hat{\eta}_0)$ be the maximizer of the null likelihood function. We assume that $\hat{\boldsymbol{\theta}}_0$ converges in probability to $\boldsymbol{\theta}_0$. Now, choose a set of functions $u_1, u_2, \ldots$ that is a basis for all functions of interest. Then use a Gram–Schmidt procedure to construct $\hat{v}_1, \ldots, \hat{v}_{n-p}$ that are linear combinations of $\gamma_1, \ldots, \gamma_p, u_1, \ldots, u_{n-p}$ satisfying (4) and (5) with $\boldsymbol{\theta}_0$ and $v_j$s replaced by $\hat{\boldsymbol{\theta}}_0$ and $\hat{v}_j$s, respectively.

For generalized linear models, the likelihood ratio statistic $\mathcal{L}_j$ can be explicitly obtained as

$$\mathcal{L}_j = 2 \sum_{i=1}^{n} [Y_i(\Delta_{ij} - \Delta_{i0}) - \{b(\Delta_{ij}) - b(\Delta_{i0})\}],$$

where, for $j = 0, \ldots, K$,

$$\Delta_{ij} = \frac{g_j(\mathbf{x}_i; \hat{\boldsymbol{\theta}}(M_j), \hat{\boldsymbol{\alpha}}_j)}{a\{\hat{\eta}(M_j)\}}.$$

Note that the maximum likelihood estimators $\hat{\boldsymbol{\theta}}(M_j)$, $\hat{\eta}(M_j)$ and $\hat{\boldsymbol{\alpha}}_j$ depend on the model used.

**3. Numerical results.** The applicability of the proposed tests is illustrated by a simulation study in a simple regression setting in Section 3.1 and by an example involving variable star data in Section 3.2. S-Plus is used for calculations.

3.1. *Simulations.* We consider normal response data

$$(6) \qquad\qquad\qquad Y_i \sim \mathcal{N}(\gamma(x_i), \eta),$$

where $x_i = (i - 1/2)/n$, $i = 1, \ldots, n$. The mean $\gamma(\cdot)$ is the parameter of interest and $\eta$ is the unknown variance parameter. In all settings the sample size was $n = 100$ and $\eta = 0.1$. We focus on testing for no effect, that is, $\gamma(x) \equiv \theta$. For the alternative models $M_j$ we take

$$\gamma_j(x) = \theta + \sum_{k \in \mathcal{K}_j} \phi_k u_k(x), \qquad j = 1, \ldots, K,$$

with $\mathcal{K}_j$ a subset of $\{1, \ldots, j\}$ and $u_k(\cdot) \equiv p_k(\cdot)$, the normalized Legendre polynomials on the interval $[1/(2n), 1 - 1/(2n)]$, $k = 1, \ldots, K$. To examine the influence of the choice of $K$, all simulations were repeated for $K = 10$ and $K = 20$.

Define $\mathrm{AIC}_j = \log L_j - m_j$ [Akaike (1974)],

$$\hat{r}_a = \underset{0 \le j \le K}{\arg\max}\, \mathrm{AIC}_j \quad \text{and} \quad \hat{r}_b = \underset{0 \le j \le K}{\arg\max}\, \mathrm{BIC}_j.$$

We compare the singleton $(B_S)$ and nested $(B_N)$ versions of the Bayes-motivated statistic $\pi_{\mathrm{BIC}}$ (Theorem 1) with some other nonparametric omnibus tests:

the tests

$$L_a = \mathcal{L}_{\hat{r}_a} \quad \text{and} \quad L_b = \mathcal{L}_{\hat{r}_b},$$

the "max-test" based on

$$M_S = \max_{0 \le j \le K} \mathcal{L}_j - 2\log K + \log\log K + \log \pi,$$

and, finally, the adaptive Neyman test $N_A$, which is based on the squared discrete Fourier transform of the residual vector from the fitted null model [Fan and Huang (2001), Section 2.1]. We also included two parametric likelihood ratio tests comparing the null model $M_0$ with the true (unknown) alternative model $M_j$ ("Oracle" test) and with the "full" model (FM) $M_K$ with $\mathcal{K}_K = \{1, \ldots, K\}$.

The tests $L_a, L_b$ and $B_N$ are all based on a sequence of nested alternative models with $\mathcal{K}_j = \{1, \ldots, j\}$, $j = 1, \ldots, K$. Score versions of $L_a$ and $L_b$ were studied in Aerts, Claeskens and Hart (2000), who established that in the present scenario $L_a$ converges in distribution to $W_{\tilde{r}}$ and $L_b$ to $W_1$, where $W_r = V_1 + \cdots + V_r$, for $r = 1, 2, \ldots, K$, $V_0 = 0$, $V_1, V_2, \ldots, V_K$ is a sequence of independent $\chi_1^2$ random variables and $\tilde{r}$ is the value of $r$ that maximizes $W_r - 2r$ over $r = 0, 1, \ldots, K$. The tests $B_S$ and $M_S$ apply singleton alternative models with $\mathcal{K}_j = \{j\}$, $j = 1, \ldots, K$. The test $M_S$ is expected to have power characteristics similar to those of $B_S$ (Theorem 5).

The definition of $\pi_{\mathrm{BIC}}$ suggests that the distributions of

$$(7) \qquad \frac{\sum_{j=1}^K \exp(V_j/2)}{1 + n^{-1/2} \sum_{j=1}^K \exp(V_j/2)}$$

and

$$(8) \qquad \frac{\sum_{j=1}^{K} n^{(1-j)/2} \exp(\sum_{i=1}^{j} V_i/2)}{1 + \sum_{j=1}^{K} n^{-j/2} \exp(\sum_{i=1}^{j} V_i/2)}$$

be used as finite sample corrected approximations to the distributions of $B_S$ and $B_N$, respectively. For $L_b$, which converges in distribution to a $\chi_1^2$ random variable, we include a corrected distribution defined as that of $W_{\tilde{r}}$, where $W_r$ is as before and $\tilde{r}$ is the value of $r$ that maximizes $W_r - r \log n$ over $r = 1, \ldots, K$.

From a simulation based on 30,000 replications, we obtained critical points (levels $\alpha = 0.01, 0.05, 0.10$) of the large sample distribution of each test statistic, except for $M_S$ and $N_A$, which are asymptotically distributed with distribution function $\exp(-\exp(-x/2))$ (Theorem 5) and $\exp(-\exp(-x))$ [Fan and Huang (2001), Theorem 1], respectively. The critical values are shown in Table 1. Although the limiting distributions of $B_N$ and $L_b$ depend on $K$, the simulations produced critical points that were identical for both values of $K$.

Table 2 shows simulated type I error probabilities for all omnibus tests based on a simulation of size 5000 with $K = 10$. The results are very similar for $K = 20$. As mentioned in Fan and Huang (2001), the approximation $\exp(-\exp(-x))$ is not so good. This was confirmed in our simulations. The simulated type I error probabilities of the adaptive Neyman test, based on the simulated critical points of Table 1 in Fan and Huang (2001), are considerably better (see last line in Table 2). The true levels of most tests are close to the nominal levels. The max-test is somewhat conservative, whereas the adaptive Neyman and the $L_a$ test are slightly liberal.

TABLE 1
*Simulated critical points of limiting null distributions*

| Test | $K$ | $\alpha = 0.10$ | $\alpha = 0.05$ | $\alpha = 0.01$ |
|------|-----|-----------------|-----------------|-----------------|
| $L_a$ | 10 | 9.393 | 13.521 | 21.028 |
| $L_a$ | 20 | 9.985 | 14.871 | 28.103 |
| $L_b$ | | 3.460 | 5.620 | 10.832 |
| $B_N$ | | 3.728 | 5.105 | 8.149 |
| $B_S$ | 10 | 8.170 | 8.724 | 9.598 |
| $B_S$ | 20 | 9.027 | 9.339 | 9.795 |

TABLE 2
*Simulated type* I *error probabilities for* $K = 10$

| Test | $\alpha = 0.10$ | $\alpha = 0.05$ | $\alpha = 0.01$ |
|------|-----------------|-----------------|-----------------|
| $L_a$ | 0.100 | 0.063 | 0.019 |
| $L_b$ | 0.102 | 0.050 | 0.010 |
| $B_N$ | 0.094 | 0.052 | 0.010 |
| $B_S$ | 0.109 | 0.055 | 0.012 |
| $M_S$ | 0.079 | 0.036 | 0.006 |
| $N_A$ | 0.125 | 0.069 | 0.017 |

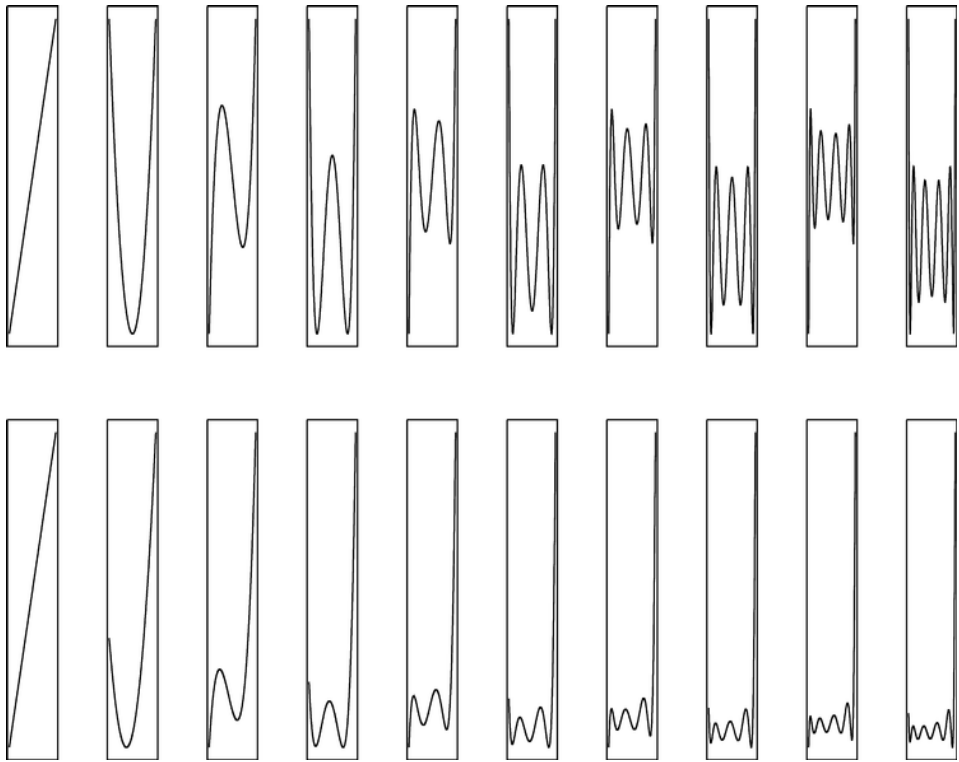

Fig. 1. *Alternative models: upper row $\gamma_m^S(x)$, lower row $\gamma_m^N(x)$, for $m = 1, \ldots, 10$.*

To examine power we consider two types of alternatives:

$$\text{(9)} \qquad \gamma_m^S(x) = u_m(x)$$

and

$$\text{(10)} \qquad \gamma_m^N(x) = \frac{1}{\sqrt{m}} \sum_{k=1}^{m} u_k(x),$$

with $m$ ranging from 1 to 10. These alternative models are ordered in the sense that they incorporate higher frequency terms as $m$ increases; for $\gamma_m^S(x)$ as single effects and for $\gamma_m^N(x)$ as nested effects (see Figure 1).

In Figure 2, power results are shown for 1000 data sets generated from the alternative models (9) and (10), respectively. In all cases, the sample size $n$ equals 100 and the level of significance is equal to 0.05. For all omnibus tests, critical points were calculated using the 5000 simulated data sets under the null hypothesis, and, hence, each omnibus test has true level very close to 0.05.

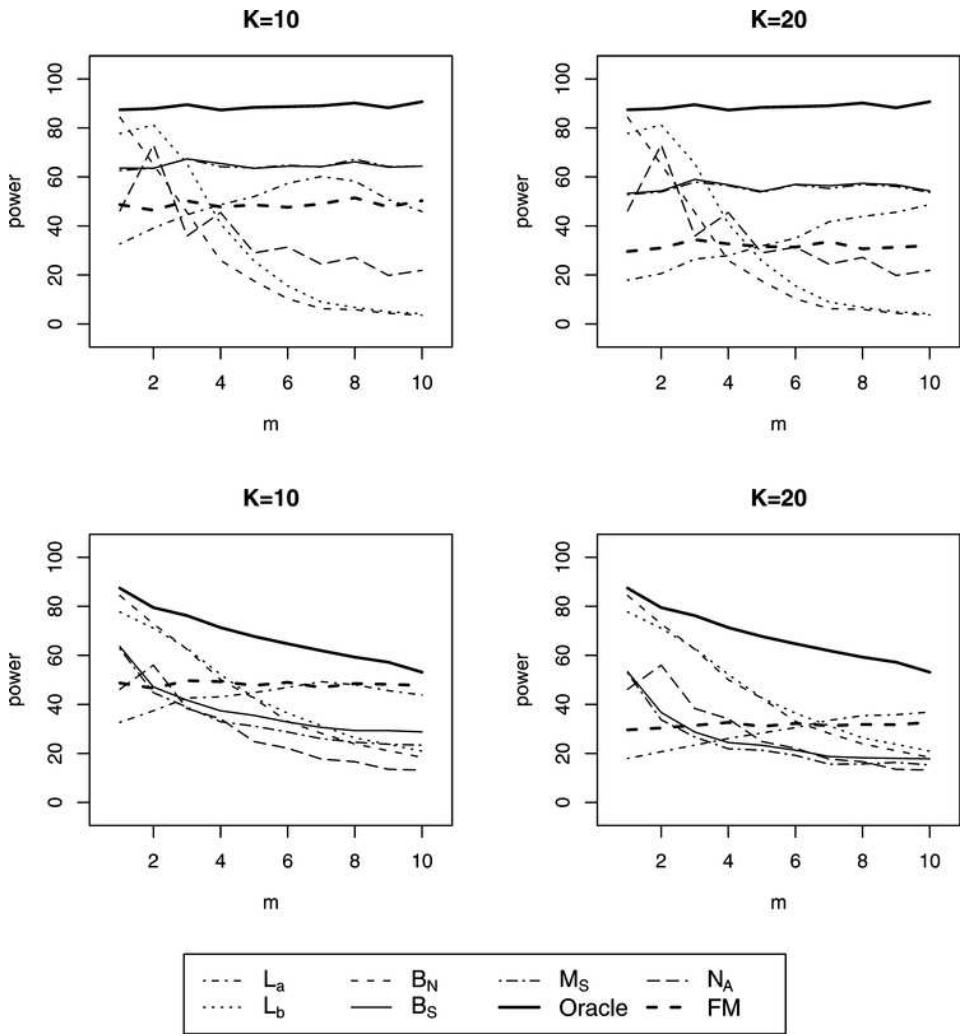

Fig. 2.  *Simulated power curves for alternative models $\gamma_m^S(x)$ (upper panels) and $\gamma_m^N(x)$ (lower panels).*

Focusing on the upper panels (single effect alternatives), four tests essentially show constant power: the Oracle test, next the singleton test $B_S$ and max-test $M_S$ with almost identical curves (as expected from Theorem 5), and the full model test. When increasing the value of $K$ from 10 to 20 (from left to right panel), the power decreases somewhat, especially for the full model test. The power characteristics of $L_b$ and $B_N$ are comparable (with some advantage for $L_b$): they have the highest power for the first lower frequency terms but their power drops down rapidly, with very comparable values for both values of $K$. The adaptive Neyman test also has a decreasing trend, but with strikingly higher powers for even alternatives. This is related to the fact that the cosine based Fourier transform terms enter the sum in the test statistic first, alternating with the sine terms. Finally, the only test with an increasing power curve is the $L_a$ test. For the single effect alternatives, the Bayesian-motivated test $B_S$ is clearly the best choice.

For the nested effect alternatives (lower panels in Figure 2), only the full model test has seemingly constant power behavior; but the higher the value of $K$ (making the test more omnibus), the less competitive this parametric approach becomes. The singleton test $B_S$ and the max-test $M_S$ are again very close and somewhat comparable to the adaptive Neyman test $N_A$. But their overall performance is rather poor. The best choices, especially for $K$ large and for alternatives $\gamma_m^N(x)$ with $m < 7$, are the Bayesian-motived test $B_N$ and the $L_b$ test. As for the single effect alternatives, the $L_a$ test seems to be a good choice for (very) high frequencies.

No single omnibus test is superior for all types of alternatives. This general statement, which is accepted as a sort of consensus by many statisticians, is confirmed by this (small) simulation study. It also shows the importance of additional knowledge, from experts in the application area, about the plausibility of certain types of alternatives.

3.2. *Analysis of data from a variable star.* Astronomers, both professional and amateur, have collected masses of data on variable stars [Mattei (1997)]. The length of time between consecutive maximum brightnesses of a star is an important quantity to astronomers since it contains information about the age and other properties of the star. We shall refer to these lengths of time as "pseudo-periods," since they tend to fluctuate substantially about the star's actual period, which is determined by fitting a periodic function to observations. Of particular interest is detecting systematic changes, or trends, in a star's period [Koen and Lombard (2001)].

Here we will apply the methodology introduced in this paper to test for period changes in the long-period variable Omicron Ceti, or Mira. Both a frequentist and a "proper" Bayesian analysis of the data will be done. The data are $(j, Y_j)$, $j = 1, \ldots, 76$, where $Y_j$ is the observed time (in days) between the $(j-1)$st and $j$th maxima on Mira's light curve. The light curve is simply Mira's brightness as a function of time. A plot of the observed pseudo-periods is given in Figure 3.

Note that we may treat $Y_1, Y_2, \ldots,$ as a time series, although the index $j$ is not actually time. A model often used by astronomers is as follows:

$$Y_j = P + \Delta_j + I_j + \varepsilon_j - \varepsilon_{j-1}, \qquad j = 1, 2, \ldots,$$

where $P$ is the mean period of the star, $\Delta_j$ is a systematic deviation from the mean period, $I_j$ represents random variation intrinsic to the star, and $\varepsilon_j$ is the error made in measuring the $j$th time of maximum brightness.

A common set of assumptions is that the $\varepsilon_j$s are i.i.d. with mean 0 and variance $\sigma_\varepsilon^2 < \infty$, the $I_j$'s are i.i.d. with mean 0 and variance $\sigma_I^2 < \infty$, and the two series are independent of each other. Our model generalizes two aspects of this one. First of all, we allow for heteroscedasticity among the $\varepsilon_j$s via the model

$$\mathrm{Var}(\varepsilon_j) = \exp(v_0 + v_1 j), \qquad j = 1, \ldots, 76.$$

This model is motivated by analysis of data from 378 variable stars by Hart, Koen and Lombard (2004), which indicates a tendency for residual variance to decrease over time, a not unexpected phenomenon since observation methods have improved with time. A second difference in our model is that we allow the $I_j$'s to follow a first order autoregressive [AR(1)] model, that is,

$$I_j = \rho I_{j-1} + Z_j, \qquad j = 2, \ldots, 76,$$

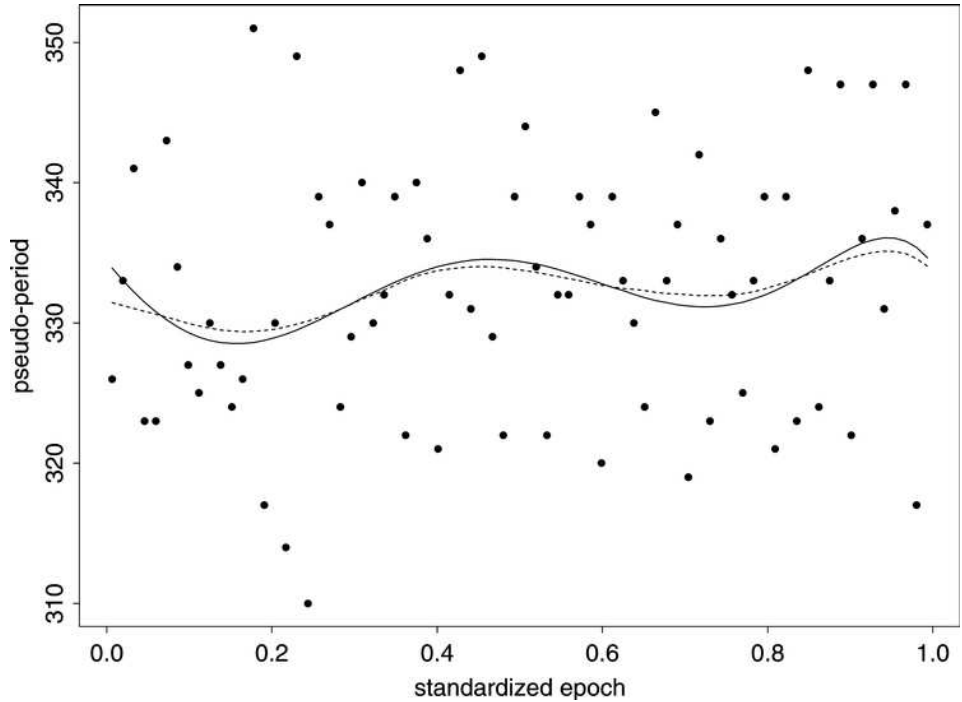

FIG. 3. *Mira pseudo-periods and two estimates of trend. The solid line is a sixth degree polynomial, the model chosen by BIC. The dashed line is a local linear smooth.*

where $|\rho| < 1$ and the $Z_j$s are i.i.d. mean 0 random variables with finite variance $\sigma_Z^2$. Our motivation for using an AR model is to circumvent a false indication of trend. It is well known that the actual size of a trend test assuming independent data is usually larger than the nominal size when the data exhibit *positive* serial correlation.

We will model the trend $\Delta_j$, $j = 1, \ldots, 76$, as a polynomial of unknown degree, and take $K = 15$ as an upper bound on the degree. To obtain a likelihood function, we assume that both the $\varepsilon_j$s and $Z_j$s are Gaussian. Therefore, our complete model says that $Y_1, \ldots, Y_{76}$ are jointly normal with means of the form

$$E(Y_j) = \beta_0 + \beta_1 j + \cdots + \beta_k j^k, \qquad j = 1, \ldots, 76,$$

and covariance matrix defined by

$$\mathrm{Cov}(Y_i, Y_j) = \begin{cases} \sigma_Z^2/(1-\rho^2) + \exp[v_0 + v_1 j] + \exp[v_0 + v_1(j-1)], & i = j, \\ \rho\sigma_Z^2/(1-\rho^2) - \exp[v_0 + v_1 \min(i,j)], & |i-j| = 1, \\ \rho^{|i-j|}\sigma_Z^2/(1-\rho^2), & |i-j| > 1. \end{cases}$$

We wish to test the hypothesis

$$H_0 : \Delta_1 = \Delta_2 = \cdots = \Delta_n = 0.$$

In our frequentist analysis, two test statistics were computed. One is $\pi_{\mathrm{BIC}}$ for the nested polynomial models with degrees $0, 1, \ldots, 15$, and the other is

$$\pi_{\mathrm{singleton}} = \left(1 + \sum_{j=1}^{15} \exp[\log(L_j/L_{j-1}) - \log(76)/2]\right)^{-1},$$

where $L_j$ is the maximized likelihood for the degree $j$ polynomial model. The components $L_j/L_{j-1}$, $j = 1, \ldots, 15$, are approximately independent of each other, with the $j$th component representing the relative increase in likelihood when stepping from a $(j-1)$st to a $j$th degree polynomial.

The values of $\pi_{\mathrm{BIC}}$ and $\pi_{\mathrm{singleton}}$ for the Mira data were 0.000161 and 0.00265, respectively. Inasmuch as these quantities are good approximations to posterior probabilities of no trend, this is already considerable evidence in favor of a trend. However, we may also use frequentist methods to judge the significance of these values. A parametric bootstrap was used to approximate the distribution of the two statistics on the assumption that $H_0$ is true. Data were generated from the estimated error model corresponding to the polynomial degree maximizing $\mathrm{BIC}_j$, $j = 0, 1, \ldots, 15$. The estimated optimal degree was 6, and the maximum likelihood estimate of $\sigma_Z^2/\exp(v_0)$ at degree 6 was 0. Essentially, this says that the experimental errors, $\varepsilon_j - \varepsilon_{j-1}$, are estimated to be so large that they completely overwhelm the intrinsic errors, $I_j$. The maximum likelihood estimate of $v_1$ at degree 6 was $-0.001816$. In our bootstrap procedure, we thus generated observations $Y_j^*$ according to

$$Y_j^* = \varepsilon_j^* - \varepsilon_{j-1}^*, \qquad j = 1, \ldots, 76,$$

where the $\varepsilon_j^*$s are i.i.d. with $\varepsilon_j^* \sim N(0, \exp(-0.001816j))$, $j = 0, 1, \ldots, 76$. (Since the distributions of our likelihood ratios are invariant to a constant mean and to $v_0$, we took these two parameters to be 0.)

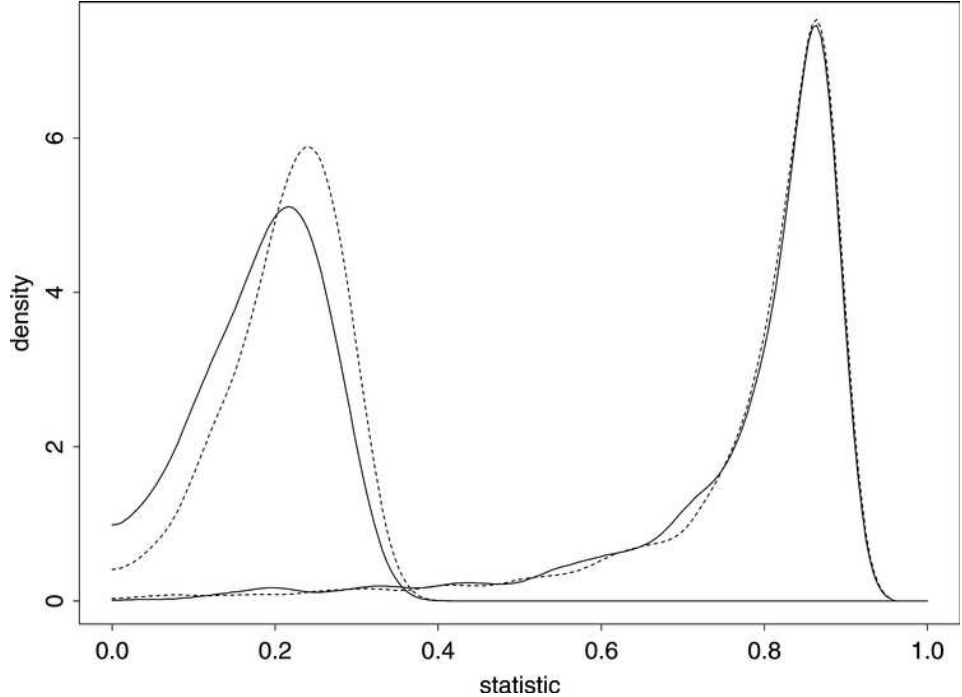

FIG. 4.  *Approximations to the distributions of $\pi_{\mathrm{BIC}}$ (right) and $\pi_{\mathrm{singleton}}$ (left). The solid lines are obtained from a Gaussian bootstrap, and the dashed lines are asymptotic distributions.*

One thousand sets of bootstrap data were generated, and on each one we computed $\pi_{\mathrm{BIC}}^*$ and $\pi_{\mathrm{singleton}}^*$ in exactly the same way that $\pi_{\mathrm{BIC}}$ and $\pi_{\mathrm{singleton}}$ were computed from the original data. Kernel density estimates for the two bootstrap distributions are shown in Figure 4. In addition, we provide estimates of the densities of $\pi_{\mathrm{BIC,asy}}$ and $\pi_{\mathrm{singleton,asy}}$, where

$$\pi_{\mathrm{BIC,asy}} = \left(1 + \sum_{j=1}^{15} \exp\left[\sum_{i=1}^{j} V_j/2 - \log(76)j/2\right]\right)^{-1},$$

$$\pi_{\mathrm{singleton,asy}} = \left(1 + \sum_{j=1}^{15} \exp[V_j/2 - \log(76)/2]\right)^{-1}$$

and $V_1, \ldots, V_{15}$ are i.i.d. $\chi_1^2$ random variables. The two latter distributions are large sample approximations to the null distributions of the two statistics.

The two approximations to the distribution of $\pi_{\text{BIC}}$ are in close agreement, while those for $\pi_{\text{singleton}}$ differ somewhat. The bootstrap distribution has a heavier left tail than the large sample approximation. Estimated $P$-values for $\pi_{\text{BIC}}$ and $\pi_{\text{singleton}}$ are 0 and $1/2000$, respectively, these being based on the two bootstrap distributions. So, the frequentist analysis provides strong evidence of a trend in the Mira pseudo-periods. Estimates of trend are seen in Figure 3.

We now describe a Bayesian analysis of the data. Priors for all model parameters were determined empirically by fitting distributions to maximum likelihood estimates for a database of 378 stars, one of which is Mira. The prior for the polynomial degree $k$ is of particular importance since the prior probability of the null hypothesis is simply the prior probability of $k = 0$. We considered three different priors for $k$: uniform over $0, 1, \ldots, 15$,

$$\pi_1(k) = \frac{1}{3.381(k+1)}, \qquad k = 0, 1, \ldots, 15,$$

and

$$\pi_2(k) = \begin{cases} 0.5, & k = 0, \\ [2(2.381)(k+1)]^{-1}, & k = 1, \ldots, 15. \end{cases}$$

The prior $\pi_1$ is a truncated version of Jeffreys' noninformative prior for an unrestricted positive integer [Jeffreys (1961), page 238], while $\pi_2$ is a modified version of $\pi_1$ that is "fair" to the null hypothesis, in that $\pi_2(0) = 0.5$. Posterior probabilities of each polynomial degree were approximated using a modification of Laplace's method that accounts for the possibility that the MLE of $\sigma_Z^2$ can occur at its lower boundary of 0. The results are given in Table 3. Regardless of which prior is used for $k$, the posterior probability of a trend is at least 0.978, and, hence, the Bayesian analysis is in basic agreement with our earlier frequentist one.

## 4. Properties of frequentist tests.
We now investigate asymptotic frequentist properties of the test statistic $1 - \pi_{\text{BIC}}$. We show how the limiting distribution of this statistic depends on the class of models $\mathcal{A}$, and we study the power of a version of the test based on singleton models (Section 2.1). It will be shown that the "singleton" test can detect local alternatives tending to the null at rate $\sqrt{\log n}/\sqrt{n}$, and that its limiting power is completely determined by the largest Fourier coefficient of the true function. Proofs of all theorems are provided in the Appendix.

Table 3

*Approximations to posterior probabilities of polynomial degrees k for Mira data. The first row is obtained using the classical BIC approximation to posterior probabilities, while the other three are based on a proper Bayesian analysis with different priors for k, as explained in the text*

| Prior | $k$ | | | | | | | | | | |
|---|---|---|---|---|---|---|---|---|---|---|---|
| for $k$ | **0** | **1** | **2** | **3** | **4** | **5** | **6** | **7** | **8** | **9** | **$\geq 10$** |
| BIC | 0.000 | 0.000 | 0.000 | 0.000 | 0.002 | 0.320 | 0.489 | 0.059 | 0.095 | 0.030 | 0.005 |
| Uniform | 0.001 | 0.003 | 0.001 | 0.000 | 0.002 | 0.107 | 0.335 | 0.109 | 0.189 | 0.129 | 0.124 |
| $\pi_1$ | 0.009 | 0.011 | 0.003 | 0.001 | 0.003 | 0.141 | 0.377 | 0.108 | 0.166 | 0.102 | 0.079 |
| $\pi_2$ | 0.022 | 0.011 | 0.003 | 0.001 | 0.003 | 0.139 | 0.372 | 0.106 | 0.163 | 0.100 | 0.080 |

4.1. *Limiting distribution under the null hypothesis.* Our first two theorems are quite general in the sense that we only make assumptions about the limiting behavior of the log-likelihood ratios $\mathcal{L}_j = 2\log(L_j/L_0)$. These assumptions hold for a great variety of likelihood models. In the sequel, $\chi^2_k$ denotes a random variable having the chi-squared distribution with $k$ degrees of freedom.

The effect of $\mathcal{A}$ is well illustrated in our first theorem, in which $\mathcal{A}$ contains finitely many models.

Theorem 1. *Let $\mathcal{A}$ be a set containing only a finite number of different models, $M_1, \ldots, M_K$, all including the null model $M_0$ as a special case. Denote by $m$ the minimal set size $m = \min_{1 \leq j \leq K}(|M_j|)$, where $|M|$ is the dimension of model $M$, and define*

$$\mathcal{K}_m = \{j \in \{1, \ldots, K\} : |M_j| = m\} = \{m(1), m(2), \ldots, m(\tilde{m})\}.$$

*We assume the following conditions hold:*

(i) *For $j = 1, \ldots, K$, the log-likelihood ratio $\mathcal{L}_j$ is bounded in probability as $n \to \infty$.*

(ii) *$(\mathcal{L}_{m(1)}, \ldots, \mathcal{L}_{m(\tilde{m})}) \xrightarrow{\mathcal{D}} (V_1, \ldots, V_{\tilde{m}})$, where $V_1, \ldots, V_{\tilde{m}}$ are jointly distributed random variables each having the $\chi^2_{m-m_0}$ distribution.*

*It then follows that*

$$n^{(m-m_0)/2}(1 - \pi_{\text{BIC}}) \xrightarrow{\mathcal{D}} \sum_{j=1}^{\tilde{m}} \exp(\tfrac{1}{2}V_j).$$

Perhaps the most important aspect of Theorem 1 is the fact that the limit distribution of $1 - \pi_{\text{BIC}}$ is completely determined by the models in $\mathcal{A}$ with the fewest parameters. In the special case where $\mathcal{A}$ is a finite sequence of nested models, Theorem 1 implies that $n^{(m-m_0)/2}(1 - \pi_{\mathrm{BIC}})$ converges in distribution to $\exp(\frac{1}{2}\chi^2_{m-m_0})$, where $m$ is the number of parameters in the smallest model in $\mathcal{A}$. This "fewest parameters" phenomenon can also be seen in the BIC-based goodness-of-fit test proposed by Ledwina (1994), and is a result of the fact that BIC consistently chooses the null model when $H_0$ is true. For more discussion on the phenomenon, see Claeskens and Hjort (2004).

Our next two theorems address cases in which the number of alternative models tends to $\infty$ with $n$. Theorem 2 is essentially a corollary to Theorem 1, and, hence, we do not provide its proof.

THEOREM 2. *Let $M_0, M_1, \ldots$ be a sequence of nested models containing numbers of parameters $m_0 < m_1 < \cdots$, respectively. Assume that under $H_0$ and as $n \to \infty$,*

$$\mathcal{L}_1 \xrightarrow{\mathcal{D}} \chi^2_{m_1 - m_0}.$$

*Furthermore, assume that, as $n \to \infty$, $\mathcal{L}_j$ is bounded in probability for each $j = 2, 3, \ldots$. Then there exists a sequence $\{K_n\}$ tending to infinity such that*

$$n^{(m_1 - m_0)/2}\left[1 - \left(1 + \sum_{j=1}^{K_n} \exp(\mathrm{BIC}_j - \mathrm{BIC}_0)\right)^{-1}\right] \xrightarrow{\mathcal{D}} \exp(\tfrac{1}{2}\chi^2_{m_1 - m_0})$$

*as $n \to \infty$.*

We now assume that the data follow a generalized linear model, as discussed in Section 2.4. We study the case where $\mathcal{A} = \mathcal{A}_K$ consists of the singleton models $M_1, \ldots, M_K$ discussed at the beginning of Section 2, and we let $K$ tend to infinity with $n$. Theorems 1 and 2 show that the asymptotic null distribution of $\pi_{\mathrm{BIC}}$ generally depends only on the models having the smallest number of elements. Therefore, our next theorem is more general than it first appears, since it also describes the limiting distribution of $\pi_{\mathrm{BIC}}$ in many cases where the alternatives consist of singletons *plus* other, larger models.

Define the statistic $S_n$ by

$$S_n = \sum_{j=1}^{K} \exp\left(\frac{n\hat{\alpha}_j^2}{2a(\hat{\eta}_0)}\right),$$

where

$$\hat{\alpha}_j = \frac{1}{n}\sum_{i=1}^{n}[Y_i - b'(g(\mathbf{x}_i; \hat{\boldsymbol{\theta}}_0))]\hat{v}_j(\mathbf{x}_i), \qquad j = 1, \ldots, K.$$

From the definition of $\pi_{\mathrm{BIC}}$, we have

$$\sqrt{n}(1 - \pi_{\mathrm{BIC}}) = \frac{\tilde{S}_n}{1 + \tilde{S}_n/\sqrt{n}},$$

where $\tilde{S}_n = \sum_{j=1}^{K} \exp(\mathcal{L}_j/2)$. The statistic $S_n$ is to $\tilde{S}_n$ as a score statistic is to a likelihood ratio statistic. The quantity $n\hat{\alpha}_j^2/a(\hat{\eta}_0)$ is known to have the same limiting distribution as the log-likelihood ratio $\mathcal{L}_j$ under the null hypothesis and general regularity conditions, which suggests that under general conditions the limiting distribution of $\sqrt{n}(1 - \pi_{\mathrm{BIC}})$ is the same as that of $S_n$. In order to simplify matters by having an explicit expression for the test statistic, we thus state Theorem 3 in terms of $S_n$.

THEOREM 3.  *Define the constants*

$$a_K = \frac{\sqrt{\pi}}{2} \cdot \frac{K}{\sqrt{\log K}} \quad and \quad b_K = \frac{Ka_K}{\sqrt{\pi}} \int_1^{\infty} \frac{\sin(x/a_K)}{x^2 \sqrt{\log x}}\, dx, \qquad K = 1, 2, \ldots.$$

*Under assumptions* A1–A8 *in the Appendix,*

$$\frac{S_n - a_K}{b_K} \overset{\mathcal{D}}{\to} S$$

*as $n$ and $K$ tend to infinity, where $S$ has the stable distribution $S_1(1, 1, 0)$, in the notation of Samorodnitsky and Taqqu (1994).*

The most interesting aspect of Theorem 3 is that the limiting distribution of $S_n$ is not normal. This results from the fact that each term $\exp(n\hat{\alpha}_j^2/[2a(\hat{\eta}_0)])$ converges in distribution to $\exp(\chi_1^2/2)$. Now, $\exp(\chi_1^2/2)$ does not have first moment finite, and, hence, the classic central limit theorem does not apply to $S_n$. However, the distribution of $\exp(\chi_1^2/2)$ *is* in the domain of attraction of the stable distribution $S_1(1, 1, 0)$, as is easily verified by checking the conditions of Theorem 1.8.1 in Samorodnitsky and Taqqu (1994).

Some remarks on the size of $K$ are in order. Ideally, we would allow $K$ to be as large as $n - p$. However, our method of proving Theorem 3 allows $K$ to be no larger than $o(n^{1/8})$. Further restrictions on $K$ may be necessary depending on the choice of basis functions. The key assumptions in this regard are A2 and A8. Suppose that $u_1, u_2, \ldots$ are trigonometric functions or Walsh functions [Golubov, Efimov and Skvortsov (1991)]. Then the bounds $B_K$ (in A2) are constant for every $K$, and no further restriction on $K$ is required. If the dimension of the covariate is 1, and $u_1, u_2, \ldots$ are Legendre polynomials, then $B_K = (\text{constant})\sqrt{K}$ [Szegö (1975), pages 68 and 184] and, again, no further restriction is needed. It is also worth mentioning that the only assumption among A1–A8 affected by the dimensionality of the covariate $\mathbf{x}$ is A2. The bounds $B_1, B_2, \ldots$ will, in some cases, tend to increase with the dimensionality of the covariate. For example, if one uses products of Legendre polynomials as basis functions, then $B_K$ will be of order $K^{d/2}$, where $d$ is the dimensionality of $\mathbf{x}$. This, of course, will further reduce the allowable size of $K$.

If the practitioner feels it necessary to choose a rather large value of $K$, and is concerned about using the large sample distribution of Theorem 3, then bootstrap methods could be used to approximate the distribution of the test statistic.

4.2. *Power against local alternatives.* Here we consider power against local alternatives, that is, alternatives that tend to the null hypothesis as $n \to \infty$. We provide rates and constants for local alternatives such that a test based on $S_n$ has power tending to 1 and another rate (and constants) such that the power tends to $p$, $\alpha < p < 1$.

THEOREM 4. *Let assumptions* A1–A8 *in the Appendix hold, and assume that the function $g$ in our generalized linear model* (GLM) *has the form*

$$g_n(x) = g(x; \boldsymbol{\theta}_0) + \left( \frac{\gamma_1 + \gamma_2 \sqrt{2 \log a_K}}{\sqrt{n}} \right) \sum_{j=1}^{m} \phi_j v_j(x),$$

*where $-\infty < \gamma_1 < \infty$ and $\gamma_2 \geq 0$ are constants. We assume that one of $|\phi_j|$ is strictly larger than all others, and define $\zeta = \max_{1 \leq j \leq m} |\phi_j| / \sqrt{a(\eta_0)}$, where $a(\eta_0)$ is the dispersion parameter in the GLM. Letting $s_\alpha$ be the $(1 - \alpha)$ quantile of the stable distribution $S_1(1, 1, 0)$ and $\Phi$ the c.d.f. of the standard normal distribution, it follows that*

$$\lim_{n \to \infty} P\left( \frac{S_n - b_K}{a_K} \geq s_\alpha \right) = \begin{cases} \alpha, & \gamma_2 < 1/\zeta, \\ \alpha + (1 - \alpha)\Phi(\gamma_1 \zeta), & \gamma_2 = 1/\zeta, \\ 1, & \gamma_2 > 1/\zeta. \end{cases}$$

It is important to note that the limiting power of the $S_n$-based test is determined by the largest Fourier coefficient of the true function. In contrast, the power of tests based only on *nested* alternatives is largely determined by the coefficients of the smallest alternative models, regardless of whether those coefficients are the largest ones. [See, e.g., Aerts, Claeskens and Hart (2000).] For this reason tests based on nested alternatives often have poor power against high frequency alternatives, since lower frequency models are the default "simplest" models. Owing to the nature of $S_n$'s null distribution, it is not too surprising that the power for $S_n$ is determined by the largest Fourier coefficient. LePage, Woodroofe and Zinn (1981) show explicitly that the limit of sums converging to a stable law is determined by the few largest summands.

The connection between $S_n$ and the largest sample Fourier coefficient becomes even clearer in the next theorem. We consider the test that rejects $H_0$ for large values of

$$(11) \qquad R_n = \max_{1 \le j \le K} \left[ \frac{n \hat{\alpha}_j^2}{a(\hat{\eta}_0)} \right]$$

and show that its limiting power against the local alternatives of Theorem 4 matches that of $S_n$. Since $R_n$ is undoubtedly more familiar to most readers than is $S_n$, this result provides a sort of benchmark for understanding the power properties of $S_n$.

THEOREM 5. *Let $R_n$ be the statistic defined in* (11), *and suppose that* A1–A8 *hold. Then if $H_0$ is true,*

$$\lim_{n, K \to \infty} P(R_n - 2 \log K + \log \log K + \log \pi \le x)$$

$$= \exp(-\exp(-x/2)) \qquad \text{for each } x.$$

*Now define $x_\alpha = -2 \log \log (1 - \alpha)^{-1}$, the $1 - \alpha$ quantile of the distribution* $\exp(-\exp(-x/2))$. *When the local alternatives of Theorem 4 hold,*

$$\lim_{n, K \to \infty} P(R_n - 2 \log K + \log \log K + \log \pi \ge x_\alpha) = \lim_{n, K \to \infty} P\left( \frac{S_n - b_K}{a_K} \ge s_\alpha \right),$$

*where the latter limit is given in Theorem 4.*

4.3. *Lindley's paradox.* Lindley's paradox refers to situations where the posterior probability of a hypothesis, $H_0$, is very high, say 0.95, and yet a frequentist test indicates strong evidence against $H_0$, in that the $P$-value for $H_0$ is small, say 0.01. Typically a frequentist does not have to deal with Lindley's paradox since he or she does not compute posterior probabilities. However, a frequentist using the tests proposed in Section 3 cannot help but notice it since the test statistic itself is a posterior probability. A level $\alpha$ test has the form

$$(12) \qquad \text{"Reject } H_0 \text{ if } \pi_n \le p_{n, \alpha},\text{"}$$

with $p_{n, \alpha} \to 1$ as $n \to \infty$. This implies that for large enough $n$, a posterior probability (for $H_0$) of, for example, 0.99, would lead to rejection of $H_0$!

For frequentists concerned with Lindley's paradox, a relevant question is "at what sample size does the paradox begin to manifest itself?" It seems reasonable to say that the paradox occurs only if we reject $H_0$ when $H_0$ is a posteriori more probable than $H_a$. Therefore, we may ask at what sample size does the critical value of test (12) become larger than $1/2$? The BIC approximation to the posterior probability of $H_0$ is

$$\pi_{\text{BIC}} = \frac{1}{1 + \sum_{k=1}^{K} n^{-(m_k - m_0)/2} L_k / L_0}.$$

Let us assume that $m_k - m_0 = 1$ for each $k$ and that $L_k/L_0$, $k = 1, \ldots, K$, are asymptotically independent, as is true for the singleton models used in Theorem 3. When $H_0$ is true, the distribution of $\pi_{\text{BIC}}$ is thus approximated by that of

$$
\text{(13)} \qquad \frac{1}{1 + n^{-1/2} \sum_{k=1}^{K} \exp(V_k/2)},
$$

where $V_1, \ldots, V_K$ are i.i.d. $\chi_1^2$ random variables. Consider the test that rejects $H_0$ at sample size $n$ and nominal level $\alpha$ when $\pi_{\text{BIC}}$ is no more than $p_{n,K,\alpha}$, the $\alpha$ quantile of the distribution of (13). For purposes of discussion, we will say that Lindley's paradox occurs when $p_{n,K,\alpha} \geq 1/2$. Of course, $\pi_{\text{BIC}}$ is only an approximation to the posterior probability of $H_0$, but Kass and Wasserman (1995) provide evidence that $\pi_{\text{BIC}}$ is an excellent approximation to $\pi_n$ for certain reference priors. To be on the safe side, we could say that $1/2 < \pi_{\text{BIC}} \leq p_{n,K,\alpha}$ is an example of Lindley's paradox in cases where one is using the appropriate reference priors.

Figure 5 displays approximations of the 95th percentile of the distribution of (13) as a function of $\sqrt{n}$ and for different values of $K$. The approximations were obtained by generating 10,000 independent values of (13). The graph indicates that for a test of nominal level 0.05, Lindley's paradox is a very large sample phenomenon when using a $K$ of 10 or more. For $K = 10$, values greater than 0.5 are not included in the rejection region until $n$ is more than 6000. On the other hand, the paradox can occur for $K = 1$ when $n$ is as small as 64. The case $K = 1$ is of particular interest since then the distribution of (13) approximates that of our statistic for testing $H_0$ against a sequence of nested alternatives.

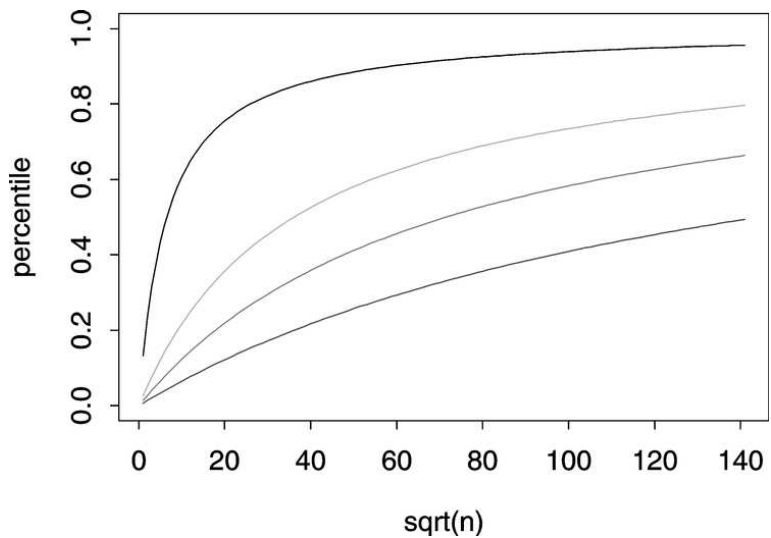

FIG. 5. *Approximate 95th percentiles of* (13). *From top to bottom, the curves correspond to* $K = 1, 5, 10, 20$, *respectively.*

A way of resolving Lindley's paradox is to use a test of the form

(14)                "Reject $H_0$ if $\pi_n \leq \min(1/2, p_{n,\alpha})$."

What effect does such a rejection region have on the power and level of the test? Typically, for all $n$ less than some $n_0$, test (14) will be identical to (12). For larger $n$, test (14) will have level of significance smaller than $\alpha$ and, indeed, tending to 0 as $n \to \infty$. Of course, the smaller rejection region will lead to an attendant reduction in power. However, in a certain sense the reduction is quite small. It can be shown that test (14) has power tending to 1 for both fixed alternatives and local alternatives tending to the null at rate $(\log n)^{\eta}/\sqrt{n}$, where $\eta > 1/2$. For local alternatives tending to 0 at rate $1/\sqrt{n}$, though, the power of (14) tends to 0. Apparently, this is a price that must be paid to avoid Lindley's paradox.

**5. Concluding remarks.** A very general means of testing the fit of a parametric function has been proposed. The parametric model is rejected if its posterior probability is too small. The test can be carried out in either a Bayesian or frequentist way. Alternatives to the null hypothesis are modeled by a sequence of models, which need not be nested. Our simulation study supports the conclusion that test validity is generally well maintained by use of an asymptotic distribution. It also shows that our proposed tests can compare favorably with other omnibus lack of fit tests.

Although some of the theory assumes the model is of the generalized linear form, the test can be used in a general likelihood context, including discrete or continuous data and multivariate data with dependence among observations. Our example using variable star data illustrates the fact that our test can accommodate dependent data. The applicability and the performance of the method in a variety of complex settings, including longitudinal and other types of clustered data, is a topic of current research. In the multiple regression case, where the covariates belong to a subset of $\mathbb{R}^d$ $(d > 1)$, a variety of sequences of alternative models, including singleton models and variations thereof, can be chosen and it is not clear which sequence is preferable or leads to optimal power characteristics. An extensive simulation study in a variety of settings can shed more light on these important practical issues. Furthermore, a score version of the proposed frequentist test can be considered, as well as robust versions of it. These variations and extensions are currently under investigation.

# APPENDIX

Following are assumptions needed in our proofs of Theorems 3 and 4:

A1. The design points $\mathbf{x}_1, \ldots, \mathbf{x}_n$ are fixed and confined to a compact subset $\mathcal{S}$ of $\mathbb{R}^d$ for all $n$.

A2. The functions $\gamma_1, \ldots, \gamma_p, u_1, u_2, \ldots$ satisfy the following assumptions:

(i) There exists $B_1^* < \infty$ such that

$$\sup_{1 \leq j \leq p, \mathbf{x} \in \mathcal{S}} |\gamma_j(\mathbf{x})| < B_1^* \quad \text{and}$$

(ii) there exists a sequence of positive constants $\{B_j : j = 1, 2, \ldots\}$ such that

$$\sup_{1 \leq j \leq K, \mathbf{x} \in \mathcal{S}} |u_j(\mathbf{x})| \leq B_K, \qquad K = 1, 2, \ldots .$$

A3. The functions $v_1, v_2, \ldots$ satisfy (4) and (5) and $\hat{v}_1, \hat{v}_2, \ldots$ are constructed from $\gamma_1, \ldots, \gamma_p, u_1, u_2, \ldots$ as described at the beginning of Section 3.

A4. Let $\mathbf{A}_{n,K}$ denote the $n \times K$ matrix with $i, j$ element $u_j(\mathbf{x}_i)$. Then we assume that the diagonal elements of $\mathbf{A}_{n,K}^T \mathbf{A}_{n,K}/n$ are all 1, and that the smallest eigenvalue of $\mathbf{A}_{n,K}^T \mathbf{A}_{n,K}/n$ is bounded away from 0 for all $n$ and $K$.

A5. The dispersion parameter $a(\eta_0)$ is positive, and the MLEs $\hat{\eta}_0$ and $\hat{\boldsymbol{\theta}}_0$ of $\eta_0$ and $\boldsymbol{\theta}_0$, respectively, are such that $E(a(\hat{\eta}_0) - a(\eta_0))^2$ and $E\|\hat{\boldsymbol{\theta}}_0 - \boldsymbol{\theta}_0\|^2$ exist and are each $O(n^{-1})$.

A6. Let $\Theta$ be the parameter space for $\boldsymbol{\theta}$. There exists a compact, connected subset $\mathcal{N}$ of $\Theta$ such that $\boldsymbol{\theta}_0 \in \mathcal{N}$ and, for each $\mathbf{x} \in \mathcal{S}$, $g(\mathbf{x}; \boldsymbol{\theta})$ is a continuous function of $\boldsymbol{\theta}$ on $\mathcal{N}$.

A7. The function $b$ is thrice differentiable with

$$\sup_{\mathbf{x} \in \mathcal{S}, \boldsymbol{\theta} \in \mathcal{N}} |b'''(g(\mathbf{x}; \boldsymbol{\theta}))| \leq B_2^*$$

for some constant $B_2^*$, the function (of $\mathbf{x}$) $b''(g(\mathbf{x}; \boldsymbol{\theta}))$ is nonnegative for each $\boldsymbol{\theta} \in \Theta$, and

$$\inf_{\mathbf{x} \in \mathcal{S}, \boldsymbol{\theta} \in \mathcal{N}} b''(g(\mathbf{x}; \boldsymbol{\theta})) > 0.$$

A8. The number of singleton models, $K$, tends to infinity with $n$ in such a way that $K \leq n^{1/8-a}$ and $B_K K^{7/2} \leq n^{1/2-a}$, where $a$ is any number such that $0 < a < 1/8$.

PROOF OF THEOREM 1. Using the explicit expression of the BIC numbers for the models under consideration, we can write the test statistic in the following form:

$$n^{(1/2)(m-m_0)}(1 - \pi_{\mathrm{BIC}}) = \frac{n^{(1/2)(m-m_0)} \sum_{j=1}^{K} n^{-(1/2)(m_j-m_0)} \exp(\mathcal{L}_j/2)}{1 + \sum_{j=1}^{K} n^{-(1/2)(m_j-m_0)} \exp(\mathcal{L}_j/2)}.$$

By definition of $m$ and assumption (i), the denominator is $1 + O_p(n^{-(1/2)(m-m_0)})$, while the numerator is equal to $n^{-(1/2)(m-m_0)} \sum_{j=1}^{\tilde{m}} \exp(\mathcal{L}_{m(j)}/2) + o_p(n^{-(1/2)(m-m_0)})$. The result now follows from assumption (ii). $\square$

PROOF OF THEOREM 3. Throughout the proof $C_1, C_2, \dots$ denote positive constants that depend on neither $n$ nor $K$. To simplify notation, we have suppressed the dependence of the $v_j$'s and $\hat{v}_j$'s on $n$.

We may express $\hat{\alpha}_j$ as

$$\hat{\alpha}_j = \tilde{\alpha}_j + e_{1j} + e_{2j},$$

where, for $j = 1, 2, \dots,$

$$\tilde{\alpha}_j = \frac{1}{n} \sum_{i=1}^{n} [Y_i - b'(g(\mathbf{x}_i; \boldsymbol{\theta}_0))] v_j(\mathbf{x}_i),$$

$$e_{1j} = \frac{1}{n} \sum_{i=1}^{n} [Y_i - b'(g(\mathbf{x}_i; \hat{\boldsymbol{\theta}}_0))][\hat{v}_j(\mathbf{x}_i) - v_j(\mathbf{x}_i)]$$

and

$$e_{2j} = \frac{1}{n} \sum_{i=1}^{n} [b'(g(\mathbf{x}_i; \boldsymbol{\theta}_0)) - b'(g(\mathbf{x}_i; \hat{\boldsymbol{\theta}}_0))] v_j(\mathbf{x}_i).$$

We may write

$$S_n = \sum_{j=1}^{K} \exp(U_{jn}/2) + \sum_{j=1}^{K} \exp(U_{jn}/2)[\exp(R_{jn}) - 1] = \sum_{j=1}^{K} \exp(U_{jn}/2) + r_{Kn},$$

where $U_{jn} = n\tilde{\alpha}_j^2/a(\eta_0)$ and $R_{jn} = (V_{jn} - U_{jn})/2$. Obviously

$$|r_{Kn}| \leq \max_{1 \leq j \leq K} |\exp(R_{jn}) - 1| \sum_{j=1}^{K} \exp(U_{jn}/2).$$

The remainder of the proof consists of two main parts:

(a) Showing that $\delta_n = \max_{1 \leq j \leq K} |\exp(R_{jn}) - 1|$ is asymptotically negligible, and

(b) obtaining the large sample distribution of $\sum_{j=1}^{K} \exp(U_{jn}/2)$.

We first consider (a). By Taylor's theorem,

$$\exp(R_{jn}) = 1 + R_{jn} \exp(\tilde{R}_{jn})$$

for $\tilde{R}_{jn}$ such that $|\tilde{R}_{jn}| \leq |R_{jn}|$, and so

$$\delta_n \leq \max_{1 \leq j \leq K} |R_{jn}| \exp(|R_{jn}|).$$

Now, for all $\varepsilon > 0$, take $\varepsilon'$ such that $\varepsilon' \exp(\varepsilon') = \varepsilon$, and so

$$(A.1) \qquad P\left(\bigcap_{j=1}^{K}\{|R_{jn}| \le \varepsilon'\}\right) \le P\left(\max_{1 \le j \le K}|R_{jn}|\exp(|R_{jn}|) \le \varepsilon\right).$$

Define

$$T_{1j} = \frac{n\tilde{\alpha}_j^2}{2}\left(\frac{1}{a(\hat{\eta}_0)} - \frac{1}{a(\eta_0)}\right), \qquad T_{2j} = \frac{n\tilde{\alpha}_j}{a(\hat{\eta}_0)}(e_{1j} + e_{2j})$$

and

$$T_{3j} = \frac{n(e_{1j} + e_{2j})^2}{2a(\hat{\eta}_0)}.$$

By (A.1) we have

$$P\left(\max_{1 \le j \le K}|R_{jn}|\exp(|R_{jn}|) > \varepsilon\right) \le \sum_{j=1}^{K}\sum_{\ell=1}^{3}P(|T_{\ell j}| > \varepsilon'/3).$$

Clearly,

$$P\left(|T_{1j}| > \frac{\varepsilon'}{3}\right) \le P\left(n^{2/3}\tilde{\alpha}_j^2 > 2\sqrt{\frac{\varepsilon'}{3}}\right) + P\left(n^{1/3}\left|\frac{1}{a(\hat{\eta}_0)} - \frac{1}{a(\eta_0)}\right| > \sqrt{\frac{\varepsilon'}{3}}\right).$$

By Markov's inequality,

$$P\left(n^{2/3}\tilde{\alpha}_j^2 > 2\sqrt{\frac{\varepsilon'}{3}}\right) \le \frac{n^{2/3}a(\eta_0)n^{-1}}{2\sqrt{\varepsilon'/3}} = \frac{\sqrt{3}a(\eta_0)}{2\sqrt{\varepsilon'}}n^{-1/3},$$

where we have used $E(Y_i) = b'(g(\mathbf{x}_i; \boldsymbol{\theta}_0))$, $\mathrm{Var}(Y_i) = a(\eta_0)b''(g(\mathbf{x}_i; \boldsymbol{\theta}_0))$ and A3. A bit of algebra shows that, for all $n$ sufficiently large,

$$P\left(n^{1/3}\left|\frac{1}{a(\hat{\eta}_0)} - \frac{1}{a(\eta_0)}\right| > \sqrt{\frac{\varepsilon'}{3}}\right) \le P\left(|a(\hat{\eta}_0) - a(\eta_0)| > \frac{a^2(\eta_0)\sqrt{\varepsilon'}}{a(\eta_0)\sqrt{\varepsilon'} + \sqrt{3}n^{1/3}}\right).$$

By Markov's inequality and A5, the last probability is $O(n^{-1/3})$. We have thus shown that $\sum_{j=1}^{K}P(|T_{1j}| > \varepsilon'/3) = O(Kn^{-1/3})$.

We turn now to the terms $T_{2j}$, for which

$$P\left(|T_{2j}| > \frac{\varepsilon'}{3}\right) \le P\left(n|\tilde{\alpha}_j||e_{1j}| > \frac{\varepsilon'a(\eta_0)}{12}\right) + P\left(n|\tilde{\alpha}_j||e_{2j}| > \frac{\varepsilon'a(\eta_0)}{12}\right)$$
$$+ P\left(a(\hat{\eta}_0) < \frac{a(\eta_0)}{2}\right).$$

The third summand in the last term is $O(n^{-1})$, independent of $j$. Letting $C_1 = \varepsilon' a(\eta_0)/12$,

$$P(n|\tilde{\alpha}_j||e_{2j}| > C_1) \leq P(n^{1/3}|\tilde{\alpha}_j| > \sqrt{C_1}) + P(n^{2/3}|e_{2j}| > \sqrt{C_1})$$
$$\leq \frac{a(\eta_0)}{C_1 n^{1/3}} + P(n^{2/3}|e_{2j}| > \sqrt{C_1}).$$

We have

$$e_{2j} = -\frac{1}{n}\sum_{i=1}^{n}\Big[(g(\mathbf{x}_i;\hat{\boldsymbol{\theta}}_0) - g(\mathbf{x}_i;\boldsymbol{\theta}_0))b''(g(\mathbf{x}_i;\boldsymbol{\theta}_0))$$
$$+ \frac{1}{2}(g(\mathbf{x}_i;\hat{\boldsymbol{\theta}}_0) - g(\mathbf{x}_i;\boldsymbol{\theta}_0))^2 b'''(\tilde{g}_i)\Big]v_j(\mathbf{x}_i),$$

where $\tilde{g}_i$ is between $g(\mathbf{x}_i;\hat{\boldsymbol{\theta}}_0)$ and $g(\mathbf{x}_i;\boldsymbol{\theta}_0)$. The orthogonality properties (4) imply that the last expression is simply $-(2n)^{-1}\sum_{i=1}^{n}(g(\mathbf{x}_i;\hat{\boldsymbol{\theta}}_0) - g(\mathbf{x}_i;\boldsymbol{\theta}_0))^2 \times b'''(\tilde{g}_i)v_j(\mathbf{x}_i)$, and so

$$|e_{2j}| \leq \left(\frac{1}{n}\sum_{i=1}^{n}(b'''(\tilde{g}_i))^2 v_j^2(\mathbf{x}_i)\right)^{1/2}\max_{1\leq i\leq n}(g(\mathbf{x}_i;\hat{\boldsymbol{\theta}}_0) - g(\mathbf{x}_i;\boldsymbol{\theta}_0))^2/2$$
$$\leq C_2\|\hat{\boldsymbol{\theta}}_0 - \boldsymbol{\theta}_0\|^2 \max_{1\leq i\leq n}\frac{|b'''(\tilde{g}_i)|}{\sqrt{b''(g(\mathbf{x}_i;\boldsymbol{\theta}_0))}}.$$

It follows that

$$P(n^{2/3}|e_{2j}| > \sqrt{C_1}) \leq P(n^{2/3}\|\hat{\boldsymbol{\theta}}_0 - \boldsymbol{\theta}_0\|^2 > C_3) + P\left(\max_{1\leq i\leq n}|b'''(\tilde{g}_i)| > C_4\right),$$

where $C_3$ and $C_4$ are defined so that $C_4$ exceeds the value $B_2^*$ in A7. We now have

$$P(n^{2/3}|e_{2j}| > \sqrt{C_1}) \leq \frac{C_5}{n^{1/3}} + P\left(\max_{1\leq i\leq n}|b'''(\tilde{g}_i)| > C_4 \cap \hat{\boldsymbol{\theta}}_0 \in \mathcal{N}\right) + P(\hat{\boldsymbol{\theta}}_0 \in \mathcal{N}^c).$$

On the event $\hat{\boldsymbol{\theta}}_0 \in \mathcal{N}$, assumption A6 implies that $\tilde{g}_i = g(\mathbf{x}_i;\boldsymbol{\theta}_{ni})$ for $\boldsymbol{\theta}_{ni} \in \mathcal{N}$, and, hence, (by A7) $P(\max_{1\leq i\leq n}|b'''(\tilde{g}_i)| > C_4 \cap \hat{\boldsymbol{\theta}}_0 \in \mathcal{N}) = 0$. Along with A5, we thus have

$$P(n^{2/3}|e_{2j}| > \sqrt{C_1}) \leq \frac{C_5}{n^{1/3}} + \frac{C_6}{n}.$$

Now consider

$$P(n|\tilde{\alpha}_j||e_{1j}| > C_1) \leq P(n^c|\tilde{\alpha}_j| > \sqrt{C_1}) + P(n^{1-c}|e_{1j}| > \sqrt{C_1})$$
$$\leq \frac{a(\eta_0)}{C_1 n^{1-2c}} + P(n^{1-c}|e_{1j}| > \sqrt{C_1}),$$

for $c$ a number in $(0, 1/2)$. We may write

$$e_{1j} = \sum_{r=1}^{p+j} (\hat{\beta}_{rj} - \beta_{rj}) \frac{1}{n} \sum_{i=1}^{n} [Y_i - b'(g(\mathbf{x}_i; \hat{\boldsymbol{\theta}}_0))] u_r(\mathbf{x}_i),$$

and so

$$e_{1j}^2 \leq \sum_{r=1}^{p+j} (\hat{\beta}_{rj} - \beta_{rj})^2 \cdot \sum_{r=1}^{p+j} \left( \frac{1}{n} \sum_{i=1}^{n} [Y_i - b'(g(\mathbf{x}_i; \hat{\boldsymbol{\theta}}_0))] u_r(\mathbf{x}_i) \right)^2$$

$$\overset{\text{def}}{=} \sum_{r=1}^{p+j} (\hat{\beta}_{rj} - \beta_{rj})^2 \cdot Z_{nj}.$$

Before proceeding, we define some matrix notation. Let $\mathbf{A}$ denote the matrix $\mathbf{A}_{n,K}$ in A4, and $\mathbf{W}$ and $\hat{\mathbf{W}}$ the $n \times n$ diagonal matrices with respective diagonal elements $b''(g(\mathbf{x}_i; \boldsymbol{\theta}_0))$ and $b''(g(\mathbf{x}_i; \hat{\boldsymbol{\theta}}_0))$, $i = 1, \ldots, n$. Matrices $\mathbf{B}$ and $\hat{\mathbf{B}}$ are the R matrices in the QR decompositions of $\mathbf{W}^{1/2} \mathbf{A}/\sqrt{n}$ and $\hat{\mathbf{W}}^{1/2} \mathbf{A}/\sqrt{n}$, respectively. We then have

$$\sum_{r=1}^{p+j} (\hat{\beta}_{rj} - \beta_{rj})^2 \leq (p+j) \max_{r,j} (\hat{\beta}_{rj} - \beta_{rj})^2 \leq (p+j) \|\mathbf{B}^{-1} - \hat{\mathbf{B}}^{-1}\|_2^2$$

$$= (p+j) \|\mathbf{B}^{-1} (\mathbf{B} - \hat{\mathbf{B}}) \hat{\mathbf{B}}^{-1}\|_2^2$$

$$\leq \frac{(p+j) \|\mathbf{B} - \hat{\mathbf{B}}\|_2^2}{\sigma^2(\mathbf{B}) \sigma^2(\hat{\mathbf{B}})},$$

where $\sigma(\mathbf{M})$ denotes the smallest singular value of matrix $\mathbf{M}$. A result of Drmač, Omladič and Veselić (1994) implies that

$$\|\mathbf{B} - \hat{\mathbf{B}}\|_2^2 \leq (8K^3 + \sqrt{2}K^2) \cdot \frac{\|\mathbf{A}^T \mathbf{W} \mathbf{A}/n\|_2 \|\mathbf{A}^T \mathbf{W} \mathbf{A} - \mathbf{A}^T \hat{\mathbf{W}} \mathbf{A}\|_2^2}{\sigma^2(\mathbf{A}^T \mathbf{W} \mathbf{A})}.$$

Assumptions A2, A4 and A7 and basic properties of matrix norms [Golub and Van Loan (1996)] now imply that

$$\|\mathbf{B} - \hat{\mathbf{B}}\|_2^2 \leq C_7 (8K^3 + \sqrt{2}K^2) \max_{1 \leq i \leq n} (b''(g(\mathbf{x}_i; \hat{\boldsymbol{\theta}}_0) - b''(g(\mathbf{x}_i; \boldsymbol{\theta}_0)))^2.$$

Combining previous results yields

$$P(n^{2(1-c)} e_{1j}^2 > C_1) \leq P\Big( n^a (p+j)(8K^3 + \sqrt{2}K^2)$$

$$\times \max_{1 \leq i \leq n} (b''(g(\mathbf{x}_i; \hat{\boldsymbol{\theta}}_0)) - b''(g(\mathbf{x}_i; \boldsymbol{\theta}_0)))^2 / \sigma^2(\hat{\mathbf{B}}) > C_8 \Big)$$

$$+ P(n^{2(1-c)-a} Z_{nj} > \sqrt{C_1})$$

$$\leq (p+j) [C_9 K^3 n^{a-1} + C_{10} n^{1-2c-a}]$$

and

$$\sum_{j=1}^{K} P\left(|T_{2j}| > \frac{\varepsilon'}{3}\right) \leq C_{11}\left(\frac{K}{n^{1/3}}\right) + C_{12}\left(\frac{K}{n^{1-2c}}\right)$$
$$+ C_{13}\left(\frac{K}{n^{(1-a)/5}}\right)^5 + C_{14}\left(\frac{K}{n^{[a-(1-2c)]/2}}\right)^2.$$

Taking $1 - 2c = (1-a)/5 = [a-(1-2c)]/2 = 1/8$ and demanding that $K = o(n^{1/8})$ ensures that the right-hand side above tends to 0.

Since $n(e_{1j} + e_{2j})^2 \leq 2n(e_{1j}^2 + e_{2j}^2)$, the term $\sum_{j=1}^{K} P(|T_{3j}| > \varepsilon'/3)$ can be bounded by a quantity that is asymptotically negligible in comparison to $\sum_{j=1}^{K} P(|T_{2j}| > \varepsilon'/3)$. Combining all the previous steps, it now follows that $\delta_n$ tends to 0 in probability as $n \to \infty$.

We turn now to step (b) in our proof. We may write

$$\frac{S_n - b_K}{a_K} = \mathcal{S}_{Kn} + \frac{r_{Kn}}{a_K},$$

where

$$\mathcal{S}_{Kn} = \frac{\sum_{j=1}^{K} \exp(U_{jn}/2) - b_K}{a_K}.$$

We will first show that $\mathcal{S}_{Kn}$ converges in distribution to a stable law, and then that $r_{Kn}/a_K$ converges in probability to 0. Now let $F_{Kn}$ be the c.d.f. of $\mathcal{S}_{Kn}$, and $F_K$ the c.d.f. of a random variable having exactly the same form as $\mathcal{S}_{Kn}$ but with $U_{1n}, \ldots, U_{Kn}$ replaced by $Z_1^2, \ldots, Z_K^2$, where $Z_1, Z_2, \ldots$ are i.i.d. random variables having the standard normal distribution. Obviously,

$$F_{Kn}(x) = F_K(x) + (F_{Kn}(x) - F_K(x)).$$

Theorem 1.8, pages 50 and 51 of Samorodnitsky and Taqqu (1994), implies that $F_K$ converges uniformly to $F$, where $F$ is the $S_1(1,1,0)$ stable law, in the notation of Samorodnitsky and Taqqu (1994).

Now $F_K$ can be written as

$$F_K(x) = P\left(\sum_{j=1}^{K} \exp(Z_j^2/2) \leq a_K x + b_K\right)$$
$$= P((Z_1, \ldots, Z_K) \in \mathcal{A}_{x,K}).$$

Likewise,

$$F_{Kn}(x) = P(\sqrt{n}(\tilde{\alpha}_1, \ldots, \tilde{\alpha}_K)/\sqrt{a(\eta_0)} \in \mathcal{A}_{x,K}).$$

Due to the convexity of the exponential function, the sets $\mathcal{A}_{x,K}$ are convex for all $x$ and $K$, and, hence, we may apply the multivariate Berry–Esséen theorem of Götze (1991) to obtain the bound

$$\sup_x |F_{Kn}(x) - F_K(x)| \leq C_{15} K \xi_{n,K},$$

for all $K \geq 6$, where

$$\xi_{n,K} = \frac{1}{(na(\eta_0))^{3/2}} \sum_{j=1}^{n} E|Y_j - b'(g(\mathbf{x}_j; \boldsymbol{\theta}_0))|^3 \left[ \sum_{i=1}^{K} v_i^2(\mathbf{x}_j) \right]^{3/2}.$$

The uniform boundedness of $b'''(g(\,\cdot\,; \boldsymbol{\theta}_0))$ (A7) now implies that

$$\begin{aligned}
\xi_{n,K} &\leq \frac{C_{16}}{n^{3/2}} \sum_{j=1}^{n} \left[ \sum_{i=1}^{K} v_i^2(\mathbf{x}_j) \right]^{3/2} \\
&\leq \frac{B_K K^{3/2} C_{16}}{n^{3/2}} \sum_{i=1}^{K} \sum_{j=1}^{n} [b''(g(\mathbf{x}_j; \boldsymbol{\theta}_0))/b''(g(\mathbf{x}_j; \boldsymbol{\theta}_0))] v_i^2(\mathbf{x}_j) \\
&\leq \frac{B_K K^{3/2} C_{17}}{n^{1/2}} \sum_{i=1}^{K} \frac{1}{n} \sum_{j=1}^{n} b''(g(\mathbf{x}_j; \boldsymbol{\theta}_0)) v_i^2(\mathbf{x}_j) \\
&= \frac{B_K K^{5/2} C_{17}}{n^{1/2}}.
\end{aligned}$$

Finally, then

$$\sup_x |F_{Kn}(x) - F_K(x)| \leq \frac{C_{18} B_K K^{7/2}}{n^{1/2}},$$

and the right-hand side of the last expression tends to 0 by A8.

Finally, consider

$$\begin{aligned}
\frac{|r_{Kn}|}{a_K} &\leq \delta_n \left[ \frac{\sum_{j=1}^{K} \exp(U_{jn}/2) - b_K + b_K}{a_K} \right] \\
&= \delta_n [O_p(1) + b_K/a_K] \\
&= o_p(1) + \delta_n b_K/a_K.
\end{aligned}$$

It is straightforward to show that $|b_K/a_K| \leq C_{19} \log K$. Examining our proof that $\delta_n$ converges in probability to 0 makes it clear that $\delta_n \log K$ does also, and, hence, the proof is complete.  $\square$

PROOF OF THEOREM 4.   For all $K > m$, we have

$$\frac{S_n - b_K}{a_K} = W_n + T_n,$$

where

$$W_n = \frac{1}{a_K} \sum_{j=1}^{m} \exp(V_{jn}/2) \quad \text{and} \quad T_n = \frac{\sum_{j=m+1}^{K} \exp(V_{jn}/2) - b_K}{a_K}.$$

Obviously

$$P(W_n + T_n \geq s_\alpha) = P(T_n \geq s_\alpha) + P(W_n + T_n \geq s_\alpha \cap T_n < s_\alpha).$$

We first consider the case where $\gamma_2^{-1} = \zeta$. Without loss of generality, suppose the largest $|\phi_j|$ is $|\phi_m|$, and consider, for any $\varepsilon > 0$,

$$P\left(\frac{1}{a_K} \sum_{j=1}^{m-1} \exp(V_{jn}/2) > \varepsilon\right) \leq P\left(\bigcup_{j=1}^{m-1} \left\{\exp(V_{jn}/2) > \frac{a_K \varepsilon}{m-1}\right\}\right)$$

$$\leq \sum_{j=1}^{m-1} P\left(\exp(V_{jn}/2) > \frac{a_K \varepsilon}{m-1}\right).$$

Now $V_{jn} = n\hat{\alpha}_j^2/a(\hat{\eta}_0)$, where

$$\hat{\alpha}_j = \frac{1}{n} \sum_{i=1}^n [Y_i - b'(g(\mathbf{x}_i; \hat{\boldsymbol{\theta}}_0))]\hat{v}_j(\mathbf{x}_i)$$

$$= \frac{1}{n} \sum_{i=1}^n [Y_i - b'(g(\mathbf{x}_i; \boldsymbol{\theta}_0))]v_j(\mathbf{x}_i) + o_p(n^{-1/2}),$$

and the last statement follows by arguing as in the proof of Theorem 3. Continuing from the last expression, we have

$$\hat{\alpha}_j = \frac{1}{n} \sum_{i=1}^n [Y_i - b'(g_n(\mathbf{x}_i))]v_j(\mathbf{x}_i)$$

$$+ \frac{1}{n} \sum_{i=1}^n [b'(g_n(\mathbf{x}_i)) - b'(g(\mathbf{x}_i; \boldsymbol{\theta}_0))]v_j(\mathbf{x}_i) + o_p(n^{-1/2})$$

$$= \frac{1}{n} \sum_{i=1}^n [Y_i - b'(g_n(\mathbf{x}_i))]v_j(\mathbf{x}_i)$$

$$+ \left(\frac{\gamma_1 + \gamma_2\sqrt{2\log a_K}}{\sqrt{n}}\right)\phi_j + O(n^{-1}\log a_K) + o_p(n^{-1/2}).$$

We may thus write

$$(\text{A.2}) \qquad \frac{\sqrt{n}\hat{\alpha}_j}{\sqrt{a(\hat{\eta}_0)}} = Z_{jn} + (\gamma_1 + \gamma_2\sqrt{2\log a_K})\frac{\phi_j}{\sqrt{a(\eta_0)}} + o_p(1),$$

where $Z_{jn}$ converges in distribution to a standard normal random variable, and we have used the fact that $\hat{\eta}_0$ is consistent for $\eta_0$ under our local alternatives. Using (A.2) and the fact that $\gamma_2|\phi_j|/\sqrt{a(\eta_0)} < 1$ for $j = 1, \ldots, m-1$, it is easy to verify that

$$P\left(\exp(V_{jn}/2) > \frac{a_K \varepsilon}{m-1}\right) \to 0$$

as $n \to \infty$ for each $j = 1, \ldots, m-1$. Combined with previous results, this implies that $\sum_{j=1}^{m-1} \exp(V_{jn}/2)/a_K$ converges to 0 in probability when $\gamma_2^{-1} = \zeta$, and, hence, the power has the same limit as

$$P(T_n \geq s_\alpha) + P(\exp(V_{mn}/2)/a_K + T_n \geq s_\alpha \cap T_n < s_\alpha).$$

Define $\tilde{T}_n$ by

$$a_K(s_\alpha - \tilde{T}_n) = \max(1, a_K(s_\alpha - T_n)).$$

Using (A.2) and some straightforward algebra, and assuming without loss of generality that $\phi_m > 0$, we have

$$Z_{mn} + (\gamma_1/\gamma_2 + \sqrt{2\log a_K}) + o_p(1)$$
$$\geq \sqrt{2\log a_K + 2\log(s_\alpha - \tilde{T}_n)} \quad \Longrightarrow \quad \exp(V_{mn}/2)/a_K + T_n \geq s_\alpha.$$

By Taylor's expansion,

$$\sqrt{2\log a_K + 2\log(s_\alpha - \tilde{T}_n)} = \sqrt{2\log a_K} + \frac{\log(s_\alpha - \tilde{T}_n)}{\sqrt{U_n}},$$

where $U_n$ is between $2\log a_K$ and $2\log a_K + 2\log(s_\alpha - \tilde{T}_n)$.

For any $0 < \varepsilon < 1$, define $I_{n,\varepsilon} = I_{(-\infty, s_\alpha - \varepsilon)}(T_n)$. Then

$$P\left(Z_{mn} + \frac{\gamma_1}{\gamma_2} + o_p(1) \geq \frac{|\log(s_\alpha - \tilde{T}_n)|I_{n,\varepsilon}}{\sqrt{2\log a_K + 2\log \varepsilon}} \cap I_{n,\varepsilon} = 1\right)$$
$$\leq P\left(Z_{mn} + \frac{\gamma_1}{\gamma_2} + o_p(1) \geq \frac{\log(s_\alpha - \tilde{T}_n)}{\sqrt{U_n}} \cap I_{n,\varepsilon} = 1\right)$$
$$\leq P\left(\exp\left(\frac{V_{mn}}{2}\right)\Big/a_K + T_n \geq s_\alpha \cap I_{n,\varepsilon} = 1\right).$$

Now, arguing as in the proof of Theorem 3, $(Z_{mn}, T_n)$ converges in distribution to $(Z, T)$, where $Z$ and $T$ are independent with standard normal and $S_1(1, 1, 0)$ distributions, respectively. Since $a_K \to \infty$, this implies that

$$\frac{|\log(s_\alpha - \tilde{T}_n)|I_{n,\varepsilon}}{\sqrt{2\log a_K + 2\log \varepsilon}} = o_p(1),$$

and so

$$P\left(Z_{mn} + \frac{\gamma_1}{\gamma_2} + o_p(1) \geq \frac{|\log(s_\alpha - \tilde{T}_n)|I_{n,\varepsilon}}{\sqrt{2\log a_K + 2\log \varepsilon}} \cap I_{n,\varepsilon} = 1\right)$$
$$\to P\left(Z \geq -\frac{\gamma_1}{\gamma_2} \cap T \leq s_\alpha - \varepsilon\right).$$

Combining previous results, and by the arbitrariness of $\varepsilon$, we now have

$$(A.3) \quad \liminf_{n \to \infty} P(\exp(V_{mn}/2)/a_K + T_n \geq s_\alpha) \geq \alpha + (1 - \alpha)\Phi(\gamma_1/\gamma_2).$$

Now, for any $\varepsilon > 0$,

$$P(\exp(V_{mn}/2)/a_K + T_n \geq s_\alpha \cap T_n < s_\alpha)$$
$$\leq P(\exp(V_{mn}/2)/a_K \geq \varepsilon \cap T_n < s_\alpha) + P(s_\alpha - \varepsilon \leq T_n < s_\alpha).$$

Arguing as we did before, the very last quantity converges to

$$(1-\alpha)\Phi(\gamma_1/\gamma_2) + P(s_\alpha - \varepsilon \leq T < s_\alpha),$$

and, hence, by the arbitrariness of $\varepsilon$,

$$\limsup_{n\to\infty} P(\exp(V_{mn}/2)/a_K + T_n \geq s_\alpha) \leq \alpha + (1-\alpha)\Phi(\gamma_1/\gamma_2).$$

Combined with (A.3), this yields Theorem 4 for the case $\gamma_2^{-1} = \zeta$.

When $\gamma_2^{-1} > \zeta$, we may show that $W_n$ tends to 0 in probability using the same argument that was applied to $W_n - \exp(V_{mn}/2)/a_K$ in the case $\gamma_2^{-1} = \zeta$.

For $\gamma_2^{-1} < \zeta$, the limiting power is at least

$$P(T_n \geq s_\alpha) + P(\exp(V_{mn}/2)/a_K + T_n \geq s_\alpha \cap T_n < s_\alpha),$$

and if we follow exactly the same steps used in the case $\gamma_2^{-1} = \zeta$, we may establish that the limiting power is 1. $\quad\square$

PROOF OF THEOREM 5.   Let $Z_1, Z_2, \ldots$ be a sequence of i.i.d. standard normal random variables, and define $d_K = 2\log K - \log\log K - \log\pi$. We first assume that $H_0$ is true. Since

$$\max_{1\leq j\leq K} U_{jn} - 2\max_{1\leq j\leq K} |R_{jn}| \leq R_n \leq \max_{1\leq j\leq K} U_{jn} + 2\max_{1\leq j\leq K} |R_{jn}|$$

and we have already shown that $\delta_n = \max_{1\leq j\leq K} |\exp(R_{jn}) - 1|$ is asymptotically negligible, it suffices to study the distribution of $\tilde{R}_n = \max_{1\leq j\leq K} U_{jn}$. Let $G_{nK}$ and $G_K$ denote the distribution functions of $\tilde{R}_n - d_K$ and $\max_{1\leq j\leq K} Z_j^2 - d_K$, respectively. The random variable $\sup_x |G_{nK}(x) - G_K(x)|$ is bounded by exactly the same quantity as was $\sup_x |F_{nK}(x) - F_K(x)|$ in the proof of Theorem 3, and, hence, we need only consider

$$P\left(\max_{1\leq j\leq K} Z_j^2 \leq x + d_K\right) = [1 - 2(1 - \Phi(\sqrt{x + d_K}))]^K$$
$$= \left[1 - \frac{2\phi(\sqrt{x + d_K})}{\sqrt{x + d_K}} + o(K^{-1})\right]^K$$
$$= [1 - e^{-x/2}K^{-1} + o(K^{-1})]^K$$
$$= \exp(-e^{-x/2}) + o(1),$$

which completes the proof in the null case.

Now assume that the local alternatives of Theorem 4 hold, and define

$$W_{1K} = \max_{1 \le j < m} U_{jn} \quad \text{and} \quad W_{2K} = \max_{m < j \le K} U_{jn}.$$

As in the proof of Theorem 4, we assume without loss of generality that the largest value of $|\phi_j|$ is at $j = m$. Three facts are key in the rest of the proof:

(i) $\tilde{R}_n$ has the same limiting distribution as $R_n$.
(ii) $P(W_{1K} \ge d_K + x_\alpha)$ converges to 0 as $n \to \infty$.
(iii) $W_{2K} - d_K$ has a limiting distribution equal to that in the null case.

Proof of (i)–(iii) is not provided here since it closely parallels arguments in the proof of Theorem 4.

Facts (i) and (ii) imply that

$$P(R_n - d_K \ge x_\alpha) = P(U_{mn} \ge d_K + x_\alpha \cup W_{2K} - d_K \ge x_\alpha) + o(1).$$

As in the proof of Theorem 4, it is easy to check that

$$\lim_{n \to \infty} P(U_{mn} \ge d_K + x_\alpha) = \Phi(\gamma_1 \zeta).$$

This along with (iii) and the fact that $U_{mn}$ and $W_{2k}$ are asymptotically independent implies that

$$\lim_{n \to \infty} P(U_{mn} \ge d_K + x_\alpha \cup W_{2K} - d_K \ge x_\alpha)$$
$$= \Phi(\gamma_1 \zeta) + \alpha - \Phi(\gamma_1 \zeta)\alpha$$
$$= \alpha + (1 - \alpha)\Phi(\gamma_1 \zeta),$$

which completes the proof. $\quad \square$

**Acknowledgments.** The authors are grateful to an Associate Editor and two referees, whose comments led to a substantial improvement in our manuscript. Finally, the authors greatly appreciate the help of Suriani Pokta in analyzing the Mira data. Her development of special Laplace approximations was a key part of our Bayesian analysis.

M. AERTS
CENTER FOR STATISTICS
LIMBURGS UNIVERSITAIR CENTRUM
UNIVERSITAIRE CAMPUS–BUILDING D
3590 DIEPENBEEK
BELGIUM
E-MAIL: marc.aerts@luc.ac.be

G. CLAESKENS
OR & BUSINESS STATISTICS
K. U. LEUVEN
NAAMSESTRAAT 69
3000 LEUVEN
BELGIUM
E-MAIL: gerda.claeskens@econ.kuleuven.be

J. D. HART
DEPARTMENT OF STATISTICS
TEXAS A&M UNIVERSITY
COLLEGE STATION, TEXAS 77843-3143
USA
E-MAIL: hart@stat.tamu.edu